\documentclass{fakeetds} 
\usepackage{amsfonts}
\usepackage{amsmath}
\usepackage{latexsym}
\usepackage{graphicx}
\usepackage{amssymb}
\usepackage{amsbsy}

\newtheorem{theorem}{Theorem}
\newtheorem{proposition}[theorem]{Proposition}
\newtheorem{lemma}[theorem]{Lemma}
\newtheorem{corollary}[theorem]{Corollary}

\def\endproof{\hfill{\vrule height4pt width6pt depth2pt}

\vspace{0.5em}

}
\def\btau{{\boldsymbol{\tau}}}
\def\ba{ {\boldsymbol{a}}}
\def\bv{ {\boldsymbol{v}}}
\def\bx{ {\boldsymbol{x}}}
\def\by{ {\boldsymbol{y}}}
\def\bu{ {\boldsymbol{u}}}
\def\bw{ {\boldsymbol{w}}}
\def\bp{ {\boldsymbol{p}}}
\def\bs{ {\boldsymbol{s}}}
\begin{document}\ETDS{1}{10}{28}{2008} 
\title{Bratteli-Vershik representations of some one-sided substitution subshifts \footnote{This research was partially supported by an NSERC grant. }}

\runningheads{ R. Yassawi}{Bratteli-Vershik representations of some one-sided substitution subshifts}
 \author{ Reem Yassawi}

\address{Dept. of Mathematics, Trent University \\
1600 West Bank Drive, Peterborough, Ontario, K9J 7B8, Canada.}

\email{ ryassawi@trentu.ca}

\recd{2009}

%\title{Representing one sided substitution shifts using adic
%systems\footnote{This research was partially supported by NSERC
%Canada.}}  

%\author{Reem Yassawi} 

%\address{Department of Mathematics, Trent
%University\\ 1600 West Bank Drive, Peterborough, Ontario, K9J
%7B8, Canada.}  

%\email{ ryassawi@trentu.ca }
%\maketitle

\begin{abstract} We study one-sided substitution subshifts, and how they can be represented using Bratteli-Vershik systems. In particular we focus on minimal recognizable substitutions such that the generated one-sided substitution 
subshift contains only one non-shift-invertible element (a {\em branch point}), and we call these substitutions {\em quasi-invertible}. We give an algorithm to check whether a substitution is quasi-invertible, 
and show that any substitution with a rational Perron value is orbit equivalent to a quasi-invertible substitution.
If the quasi-invertible substitution is left proper, then its subshift is equal to a substitution subshift where the original branch point is the new substitution fixed point. We use these results to prove that any reasonable quasi-invertible substitution subshift has a Bratteli-Vershik representation. We also give an example of a pair of substitutions whose 2-sided subshifts are topologically conjugate, while their 1-sided subshifts are not.

\end{abstract}

\section{Introduction}
Relating Bratteli-Vershik, or adic, systems whose Bratteli diagrams
have specific properties, to families of minimal 2-sided subshifts,
started in \cite{f} and \cite{dhs}, where stationary proper adic
systems were shown to be the union of the class of (two-sided)
primitive aperiodic substitution subshifts, and the class of
stationary odometers. Since then, proper adic systems whose Bratteli
diagram have the equal path number property were shown to be Toeplitz
systems in \cite{gj1}, and the adic systems realising interval exchange
maps, linearly recurrent systems, and non-minimal substitution
systems, are also described in \cite{gj2}, \cite{cdhm} and \cite{bkm}
respectively.  Implicitly used was the fact, in \cite{hps}, that any
 minimal homeomorphism of a Cantor
space is topologically conjugate to a proper  adic system. 

 If $T$ is assumed merely to be continuous, and minimal, there does not 
seem to be
straightforward realization of $(X,T)$ as an adic system, and some
constraints on the adic representation have to be imposed. For example, while 
Bratteli diagrams
which are no longer proper have to be considered, minimality of $T$ 
implies that there have to be at least as many maximal elements in 
$X_{\cal B}$ as minimal elements.

In this article we address the question of how these results can be
generalized to representing one-sided (non-$\sigma$-invertible)
substitution subshifts.  Conjugacy of a pair of two-sided substitution subshifts does not imply conjugacy of the corresponding one-sided subshifts (see Example \ref{notonesidedconjugate}), so it may not be the case that an adic representation for the one-sided subshift can be obtained directly from the adic representation of the two-sided subshift. We deal with the special case where the
substitution subshift $(X_{\btau},\sigma)$ has a unique non-invertible
point, and call these subshifts {\em quasi-invertible}. Subshifts with
this property have been studied: Sturmian sequences, and more
generally, {\em Arnoux-Rauzy sequences} (\cite{AR}) generate
quasi-invertible subshifts - the condition (*) in \cite{AR} ensures
this, and in many cases these sequences are fixed points of
substitutions (see \cite {cmps} for a characterization of when
Sturmian subshifts are substitution subshifts; if the `S-adic'
expansion of a sequence in \cite{AR} is periodic, then it is the fixed
point of a substitution, for example the tribonacci substitution). All
these subshifts have nice geometric representations as toral
rotations, or interval exchange maps.  Many (but not all) Pisot
substitutions that have been studied in the literature are
quasi-invertible. Any substitution with a rational Perron eigenvalue is orbit equivalent to some quasi-invertible substitution (see Example \ref{orbitequivalence}).
It would be interesting to study whether
quasi-invertibility imposes any spectral constraints. While not all
have discrete spectrum (see Example \ref{chacon}: the Chacon substitution is quasi-invertible) spectral
multiplicity may be bounded.

 As invertible substitution subshifts were considered in \cite{dhs}, bilateral recognizability (first introduced in \cite{Mo1}) was a sufficient assumption, and this is a property that all primitive aperiodic substitutions enjoy. Here we have to  restrict ourselves to unilaterally recognizable substitutions. Thus we do not realize  adic systems, that arise from a non-recognizable substitution, as substitution subshifts. As we want to obtain adic representations
of some topological rank one systems arising from a substitution, we work with 
 minimal (not primitive)
substitutions. We also assume that all powers of the substitution $\btau$ are injective. We'll call substitutions which satisfy these properties {\em reasonable}.  We represent these systems with a stationary Bratteli diagram which has a unique maximal
element.
We do this by generalising the techniques in \cite{dhs}, so that a reasonable quasi-invertible
substitution subshift is conjugate either to an adic system, or a `pinched' adic system.

 In Section \ref{preliminaries} we
 give definitions, and some background of required results for substitution subshifts.   In
Section \ref{suffix}, we characterize quasi-invertible substitutions (Lemma \ref{nonfixedlemma1a} and the remarks following it), and show that any left proper quasi-invertible substitution subshift is equal to a substitution subshift where the branch point is fixed by the substitution (Corollary \ref{nonfixedcorollary}. We believe that this result is still true for non-left proper subshifts -see Example \ref{adic}).
In Section
\ref{Bratteli}, we set up background and notation for adic systems, and use the results of Section \ref{suffix} to obtain adic representations for left proper quasi-invertible 
substitutions (Theorem \ref{quasiinvertibleleftproperareadic}). We extend this result to arbitrary quasi-invertible substitutions (Theorem \ref{quasiinvertiblearealmostadic}) in Section \ref{inducedsubstitutions}. 

\ack{The author thanks Ian Putnam for help with Theorem \ref{quasiinvertiblearealmostadic}.}

\section{Preliminaries \label{preliminaries}}
\subsection{Notation} 

 Let ${\cal A}$ be a finite alphabet; elements in ${\cal A}$ will be
 denoted by $a,\, b,\, c\,$ etc.  A {\em word} $\bw$ (denoted using boldface) from ${\cal
 A}$ is a finite concatenation of elements from ${\cal A}$.  Let
 ${\cal A}^{*}$ denote the set of words (including the empty word) from ${\cal A}$, and ${\cal A}^{+}$ the set of nonempty words from ${\cal A}$.   If ${\bw}=w_{1}\ldots w_{m}$ and ${\bw^{*}}=w_{1}^{*}\ldots w_{n}^{*}$,
 then ${\bw\bw^{*}}:= w_{1} \ldots w_{m}w_{1}^{*}\ldots
 w_{n}^{*}$. Define $\bw_{[i,j]}= w_{i}w_{i+1}\ldots w_{j}$.  If
 $\bw= w_{1}\ldots w_{m}$, the {\em length} of $\bw$, written  $|{\bw}|$, is $m$. Given words $\bw$, ${\bw^{*}}$, say that $\bw$
 {\em occurs} in ${\bw^{*}}$, or  $\bw$ {\em is a subword of} ${\bw^{*}}$,
if $\bw= {\bw^{*}}_{[i, i+|{\bw}|-1]}$ for some $i$.  If $|\bw|=m$, define the {\em
 k-shift} $\sigma^{k}(\bw)=\bw_{[k+1, m]} $. We say that ${\bw^{*}}$ is
 a {\em prefix} of $\bw$ if $\bw_{[1,|{\bw^{*}}|]}= {\bw}^{*}$, similarly ${\bw^{*}}$ is a {\em suffix} of $\bw$ if
 $\bw_{[|\bw|-|{\bw^{*}}|+1,|\bw|]}= \bw^{*}$. When we speak of `the' prefix (suffix) of a word $\bw$, we will mean the first (last) letter of $\bw$.

 Let ${\mathbb N}:=\{0,1, \ldots\}$. If ${\mathbb M} = {\mathbb Z}$ or
 ${\mathbb N}$, then the space of all ${\mathbb M}$-indexed sequences
 from ${\cal A}$ is written as ${\cal A}^{\mathbb M}$, and a
 configuration $\bx\,\in {\cal A}^{\mathbb M}$ is written $\bx=
 (x_{m})_{m\, \in {\mathbb M}}$. Let ${\cal A}$ be endowed with the
 discrete topology and ${\cal A}^{\mathbb M}$ with the product
 topology (this topology is also generated by the Hamming metric),
 then ${\cal A}^{\mathbb M}$ is a {\em Cantor space}: a
 zero-dimensional perfect compact metric space.  If $b\,\in {\cal A}$
 and $j\, \in \mathbb M$, the clopen sets $\{\bx: \,\,x_{j}=b\}$ form
 a countable basis for the topology on ${\cal A}^{\mathbb M}$.  The
 (left) {\em shift map} $\sigma: {\cal A}^{\mathbb M} \rightarrow
 {\cal A}^{\mathbb M}$ is the map defined as $(\sigma(\bx))_{m} =
 x_{m+1}$.  If $n\geq 2$ and ${\bf w}=w_{1}w_{2}\ldots w_{n}$, we
 define $\sigma({\bf w})=w_{2}\ldots w_{n}$. We say that $X$ is a {\em
 subshift of $({\cal A}^{\mathbb M},\sigma)$} if it is a closed
 $\sigma$-invariant subset of ${\cal A}^{\mathbb M}$. A {\em language}
 is a non-empty collection of words from ${\cal A}^{*}$.  If $\bx\,
 \in {\cal A}^{\mathbb M}$, the word $\bw$ {\em occurs in} $\bx$, or
 {\em is a word in} $\bx$ if there is some $m\, \in {\mathbb M}$ such
 that $x_{m+1} \ldots x_{m+|\bw|} = \bw$. Let ${\cal L}_{\bx}$, the
 {\em language of } $\bx$ be the collection of words occurring in
 $\bx$. If ${\cal L}$ is a collection of words closed under the
 operation of taking subwords, then the set $X^{\mathbb M}_{\cal L}$
 associated with ${\cal L}$ is the set of   sequences in ${\cal
 A}^{\mathbb M}$ all of whose words belong to ${\cal L}$. If ${\cal
 L}= {\cal L}_{\bx}$, we write $X^{\mathbb M}_{\bx}$, dropping the
 `$\mathbb M$' when the context makes it clear, or when a statement is
 true either way.

\subsection{Substitutions}

A {\em substitution} is a map $\btau: {\cal A} \rightarrow {\cal
A}^{+}$.
We extend $\btau$ to a map $\btau: {\cal
A}^{+}\rightarrow {\cal A}^{+}$ by concatenation: if ${\ba}=
a_{1}\ldots a_{k}$, then $\btau({\bf a}):= \btau(a_{1})
\ldots \btau(a_{k})$. In this way iteration $\btau^{n}$ is well defined. 
The substitution $\btau$ is extended to a map $\btau: {\cal A}^{\mathbb
N} \rightarrow {\cal A}^{\mathbb N} $ defined by $ \btau (\bx) := 
\btau (x_{0}) \,\btau (x_{1})\ldots $, and also  $\btau: {\cal A}^{\mathbb
Z} \rightarrow {\cal A}^{\mathbb Z} $ defined by $ \btau (\bx) :=
\ldots \btau (x_{-1}) \cdot \btau (x_{0}) \,\btau (x_{1})\ldots $. 
 We say $\btau$ is {\em left (right) proper} if there exists $l$ in ${\cal A}$
such that $l$ is a prefix (suffix) of $\btau (a)$ for each $a\,\in {\cal A}$.
 If $\btau$  is both left and right
proper, it is called {\em proper}. We say ${\btau}$ {\em suffix-permutative} if the set of the suffixes 
 of $\{\btau(a)\}_{a\in {\cal A}}$ is a permutation of ${\cal
A}$. Also $\btau$ is {\em injective} if $\btau(a)\neq \btau(b)$ whenever $a\neq b$.

We say ${\btau}$ is {\em primitive} if there exists a
positive integer $k$ such that for any $a\, \in {\cal A}$, all letters
of ${\cal A}$ appear in $\btau^{k}(a)$ (this requires that for some
letter $a$, $|\btau(a)|>1$). If ${\cal L}_{\btau}$ is the language generated by the words
$\{\btau^{n}(a): n\, \in {\mathbb N}^{+},\,\,a
\, \in {\cal A}\}$, then  $X^{\mathbb M}_{\btau}:= 
X^{\mathbb M}_{{\cal L}_{\btau}}$ is closed and $\sigma$-invariant, so that  
$(X^{\mathbb M}_{\btau}, \sigma)$ is a  (substitution) subshift, with the subspace topology.   Henceforth we will be
working mainly with one sided substitution subshifts, and unless otherwise indicated, will assume this, and write $(X_{\btau},\sigma)$ instead of $(X^{\mathbb N}_{\btau},\sigma)$. A {\em fixed point} of $\btau$ is a sequence ${\bu
}
\, \in {\cal A}^{\mathbb N}$ such that $\btau({\bu}) = {\bu}$. If ${\bu}\, \in {\cal A}^{\mathbb N}$ is a fixed point for $\btau$, then
$\btau(u _{0})$ starts with $u_{0}$. Conversely, if there exists a
letter $l$ such that $l$ is a prefix of $\btau(l)$, then ${\bu}:=\lim_{n\rightarrow \infty} \btau^{ n}(l)$ is the unique fixed
point satisfying $u_{0}=l$.  If $\btau$ is primitive, then the subshift generated by $\bu$, i.e. $\overline{\{\sigma^{n}({\bu})\}}$, equals $ X_{\btau}$.
Using the pigeonhole principle, there exists $n\geq 1$ such that
$\btau^{ n}$ has at least one fixed point. Though we do not  assume that $\btau$ is primitive, we will 
assume that $\btau$ has a fixed point $\bu$ with ${\cal L}_{\btau} =
{\cal L}_{\bf u}$, and $|\btau^{n}(u_{0})|\rightarrow \infty$ as $n\rightarrow\infty$.
 When such a
fixed point exists, which we henceforth assume, we will call it a {\em
generating} fixed point. Note that if $\btau$ is left proper and has a generating fixed point, then it is primitive.

The continuous mapping $T:X\rightarrow X$ of the compact metric space
$X$ is {\em minimal} if for each $x$ in $X$, $\{T^{n}(x): n \geq 0\}$
is dense in $X$. Primitive substitutions are minimal.  A sequence $\bx\, \in {\cal A}^{\mathbb M}$ is {\em
almost periodic} if for any neighbourhood $V$ of $\bx$, the set $\{ n:
\sigma^{n}(\bx) \,\in V\}$ has bounded gaps in ${\mathbb M}$. The generating fixed point  $\bu$
is almost periodic if and only if  $(X_{\btau}, \sigma)$ is minimal. (see
\cite{Qu}, Theorem 4.12). 
In what follows we will assume that our substitution $\btau$ has a
generating fixed point $\bu$, and 
that the resulting subshift $(X_{\btau},\sigma)$ is
minimal. We call such substitutions {\em minimal}.
Minimal substitutions enjoy many established properties that primitive
substitutions do, and minor modifications of proofs of facts for
primitive substitutions can be made to ensure that the same results hold
for minimal substitutions. In particular unilateral recognizability,
(see below) is defined for primitive substitutions. The necessary and sufficient
condition that exists for a primitive substitution to be unilaterally recognizable can be modified with the following two lemmas (which appear in \cite{Mo1} for primitive substitutions):

\begin{lemma}
\label{minimal2}
Suppose that $\btau$ is minimal with generating fixed point ${\bf u}$.
Fix a natural number $l$. If $\bw$ occurs in $\bu$, $|\bw|=l$, and
$\btau^{p}(\bw)=B_{p}\btau^{p}(\bw^{*})B_{p}'$, with $\bw^{*}$ occurring in $\bu$ and
$|\bw^{*}|=l^{*}$, then there is  $C$ such that $l^{*}\leq C\, l$, where $C$ is independent of $p$.
\end{lemma}
{\bf Proof:}  Suppose not. Then for all $N$ there is a $p\,\in \mathbb N$, a word 
${\bw}^{*}\, \in {\cal L}_{\btau} $ such that $\btau^{p}(\bw)=B_{p}\btau^{p}(\bw^{*})B_{p}'$
with 
$|{\bw^{*}}|=l^{*}$ and $l^{*}>Nl$. If $\bw=w_{1}\ldots w_{l}$ and 
$\bw^{*}=w_{1}^{*}\ldots w_{l^{*}}^{*}$, then 
\[\left| \prod_{i=1}^{l} \btau^{p}(w_{i} )\right|\geq
\left| \prod_{i=1}^{l^{*}} \btau^{p}(w_{i}^{*}) \right|\, .\] Since
$\btau$ is minimal, $\bu$ is almost periodic, which means that all
letters occur in $\bu$ with bounded periodicity. Hence if for  $\alpha \in {\cal A}$  we let
$n(\alpha)$ denote the number of occurrences of $\alpha$ in $\bw$, then
there are more than  $n(\alpha)$ occurrences of $\alpha$ in $\bw^{*}$ for each $\alpha$,
 if $N$ is large enough. This contradicts the inequality above.
\endproof

\begin{lemma}
\label{minimal1}Suppose that $\btau$ is minimal with generating fixed point $\bu$.
Fix a natural $k$. The sequence 
\[L_{p}:=
\frac{\max \{|\btau^{p}(\bw)|: \bw\, \in {\cal L}_{\btau}, |\bw|=k\}}
{\min\{|\btau^{p-1}(a)|:a\, \in {\cal A} \mbox{ and }\lim_{n\rightarrow\infty}|\btau^{n}(a)| = \infty\}}\] is bounded.
\end{lemma}

{\bf Proof:} 
Let ${\cal A}_{\infty}:=\{ a \in {\cal A}: \lim_{n\rightarrow\infty}
|\btau^{n}(a)| = \infty\}= \{\alpha_{1}, 
\ldots \alpha_{|{\cal A}_{\infty}|} \}.
$
 It is sufficient to show that for some $C$, 
$|\btau^{p}(\alpha)|\leq C |\btau^{p}(\beta)|$ for all natural $p$ and all $\alpha, \beta \in {\cal A}_{\infty}$. If $M$ is the 
$|{\cal A}_{\infty}|$-sized square matrix whose $i,j$-indexed entry is the number of occurrences of $\alpha_{j}$ in $\btau (\alpha_{i})$, then without loss of generality, $M$ can be assumed to have all positive entries, and the Perron-Frobenius theorem  \cite{Se} implies that 
$|\{ \beta  \in \btau^{p}(\alpha): \beta \in  {\cal A}_{\infty}\}| =C_{\alpha}\lambda^{p}+o(\lambda^{p})$, where $C_{\alpha}$ is strictly positive and $\lambda$ is the dominant eigenvalue of $M$. Thus there is a $C$ such that 
\[|\{ \gamma \in \btau^{p}(\alpha): \gamma \in  {\cal A}_{\infty}\}|
\leq C
|\{ \gamma  \in \btau^{p}(\beta): \gamma \in  {\cal A}_{\infty}\}|
\] whenever $\alpha, \beta\in {\cal A}_{\infty}$. The result follows from the almost periodicity of $\bu$.
\endproof

{\em Unilateral recognizability} was a condition introduced in \cite{Ho} and \cite{Qu} which ensured that one could find clopen generating partitions for $(X_{\btau},\sigma)$. Specifically, the set of (1)-cuttings of $\bu$ are
\[E:=\{0\}\cup\{ |\btau (\bu_{[0,p-1]})|: p>0\}\,.\]
We then say that $\btau$ is {\em (unilaterally) recognizable} if there
exists some $L$ so that if $\bu_{[i,i+L-1]}= \bu_{[j,j+L-1]}$ and
$i\, \in E$, then $j\, \in E$. A word $\bw$ occurring at $i$ and $j$ in ${\bf
u}$  has the {\em same 1-cutting} at $i$ and $j$,
with $i<j$, if $E\cap \{i,\ldots i+|\bw|\} = j-i + E\cap \{ j ,\ldots j+|\bw|\}.$

In Theorem 3.1, \cite{Mo1}, primitive recognizable substitutions are characterized, and Lemmas \ref{minimal2} and \ref{minimal1} can be used to extend this theorem to minimal  recognizable substitutions:

\begin{theorem}
\label{mosserecognizability}
Let $\btau$ be minimal  with generating  fixed point $\bu$. $\btau$ is not
recognizable if and only if for each $L$, there exists a word $\bw$ of
length L, and two elements $a, \, b \,\in {\cal A}$ such that 
\begin{itemize}
\item $\btau(b)$ is a proper suffix of $\btau(a)$.
\item The words $\btau(a) \bw$ and $\btau (b) \bw$ appear in 
$\bu$ with the same 1-cutting of $\bw$.\endproof
\end{itemize}
 \end{theorem}

 In particular,
suffix permutative substitutions are recognizable. It is
straightforward to show that if $\btau$ is recognizable, and injective,  then 
$\btau^{n}$ is recognizable for each $n$.

In Lemmas 2 and 3, \cite{Ho}, the following is shown for primitive
recognizable substitutions, and with the modifications above, the proof extends to minimal substitutions. Note that $\btau (X_{\btau}) \subset X_{\btau}$.

\begin{proposition}\label{hosttopologyprop}

 Let $\btau$ be a minimal,  recognizable  substitution with 
generating fixed point $\bu$. 

\begin{enumerate}
\item
$\by\, \in \btau (X_{\btau})$ if and only if $\sigma^{n_i} \bu
\rightarrow \by$ where $n_{i}\, \in E$ for all large $i$.
\item
$\btau (X_{\btau})$ is clopen in $X_{\btau}$.
\item
$\sigma^{p}(\btau (x)) \ \in \btau(X_{\btau})$ if and only if $p=|\btau(\bx_{[0,r]})|$ for some $r$ so that $\sigma^{p}(\btau(x)) = 
\btau(\sigma^{r+1}(\bx))$.
\end{enumerate}
If $\btau$ is also injective, then
\begin{enumerate}
\setcounter{enumi}{3}
\item
$\btau:X_{\btau}\rightarrow \btau(X_{\btau}$) is a homeomorphism. \endproof
\end{enumerate}
\end{proposition}

We remark here that the previous proposition also implies that for each $n\geq 1$, $\btau^{n}:X_{\btau}\rightarrow {\btau}^{n}(X_{\btau})$ is also a homeomorphism, though as a substitution it is not necessarily injective. 

Let us call substitutions which are recognizable, minimal, and all of whose powers are injective {\em reasonable}. In some of what follows (excepting Corollary \ref{clopenpartition}, and all results that require it), we need only assume that $\btau$ is injective, but since our final results assume that all powers of $\btau$ are injective, we will make this blanket assumption everywhere.

\begin{corollary}
\label{almostuniquerep}
Suppose that $\btau$ is reasonable.
Every $\by\,\in \btau (X_{\btau})$ can be written in a unique way as ${\bf
y}= \btau (\bx)$ with $\bx\,\in X_{\btau}$. Every $\by \, \in X_{\btau}
\backslash \btau (X_{\btau})$ can be written as $\by= \sigma^{k} \btau
(\bx)$ where $\sigma(\bx)$ is unique, and $0<k<|\btau (x_{0})|$.

\end{corollary}

{\bf Proof:} The first assertion follows immediately from Proposition
\ref{hosttopologyprop}. If  $\by \,\in X_{\btau}\backslash \btau (X_{\btau})$ then by 
Part 1 of Proposition \ref{hosttopologyprop}, there is a sequence
$(n_{i})$ of integers not in E, with $\sigma^{n_i}(\bu) \rightarrow
\by$.  Now each $n_{i}= k_{i} + r_{i}$ where $k_{i} \, \in E$, 
and the positive 
$r_{i}$ is strictly less than the successor of $k_{i}$ in $E$. Since
$E$ is almost periodic, by dropping to a subsequence if necessary, we
assume that there is some fixed positive $r$ such that each 
$n_{i}= k_{i} +r$. Now we can assume that $\sigma^{k_i}(\bu)$ converges to
some $\bx \, \in \btau(X_{\btau})$. So $\by = \sigma^{r}\bx$.

Now suppose that $\by = \sigma^{r} \bx = \sigma^{r'}\bx'$, where  
if $\bx = x_{0}x_{1}\ldots$ then $r< |\btau(x_{0})|$, and similarly for 
$r'$. Let $l= |\btau (x_{0})|-r$, and define $l'$ similarly.

If $l=l'$, then $\sigma^{r+l} \btau(\bx)= \btau (\sigma\bx)\, \in
\btau (X_{\btau})$ and similarly $\sigma^{r'+l} \btau (x') = \btau (\sigma
\bx') \, \in \btau (X_{\btau})$, so $\btau (\sigma \bx) = \btau (\sigma
\bx')$. Thus $\sigma \bx = \sigma \bx'$.

Finally we will show that $l\neq l'$ leads to a contradiction. Suppose
$l<l'$. Then $\sigma^{r+l} \btau(\bx)= \btau (\sigma\bx)\, \in \btau
(X_{\btau})$, and so $ \sigma^{r'+l} \btau(\bx') = \sigma^{r+l} \btau(\bx)\in 
\btau (X_{\btau})$. Thus by part 3 of Proposition \ref{hosttopologyprop}, $r'+l = 
|\btau (\bx'_{[0,j]})|$ for some $j$. This contradicts  
$r'+l< r'+l' = |\btau (x'_{0})|$.
\endproof

\section{Quasi-invertible substitutions \label{suffix}}

%The following proposition classifies which points have several
%$\sigma$-preimages in $(X,\sigma)$.

We investigate the number of $\sigma$-preimages that a
$\btau$-fixed point $\bu$ can have. Say that $a_{1}$ is {\em part of a
suffix cycle} $a_{1}, a_{2},\ldots a_{m}$ if
$\btau(a_{i})$ has $a_{i+1}$ as a suffix for $1\leq i<m$, and $\btau (a_{m})$ has $a_{1}$ as a suffix.

\begin{proposition}
\label{inverseimagesofu}
Suppose that $\btau$ is a reasonable
substitution with  generating fixed point $\bu = u_{0}u_{1}u_{2}\ldots$.
 Then $a\bu\, \in X_{\btau}$ if and only if  $a$ is part of a suffix cycle $a_{1},a_{2},\ldots a_{m}$, and $a_{i}u_{0}\, \in {\cal L}_{\btau}$ for some $i\, \in \{1,\ldots m\}$.
\end{proposition}

{\bf Proof:} If for some $i$, $a_{i}u_{0}\, \in {\cal L}_{\btau}$, then  since
$\bu$ is generating, there exists an $n$ such that $a_{i}u_{0}\,\in
{\btau}^{n}(u_{0})$, so $a_{i+1}\,\btau(u_{0})\, \in \btau^{n+1}(u_{0})$, 
and $a_{i}\btau^{nm}(u_{0})\,\in {\btau}^{nm}(u_{0})$ for each $n\geq 1$. It follows that  $a_{i}\bu\, \in X_{\btau}$, and that this is true for each $a_{j}\, \in \{a_{1},\ldots a_{m}\}.$% follows that
%= {\bf p}b{\bf s}$ for some words ${\bf p}$ and $ {\bf
%s}$. Thus if  $\btau(b)$ has  $a$ as a suffix,
%% where $\lim_{n\rightarrow \infty} |{\btau}^{n}({\bf s})|=
%\infty$. 
%then $a\, l \, \in {\btau}^{n+1}(l)\,\in {\cal L}_{\btau}$,  and  if 
% $a$ is part of a suffix cycle $a,b_{1}, \ldots , b_{k-1}$, then
% $b_{i}\, \l\, \in X$ for each  $i$.
%Since $b_{k-1}\,l\,\in {\cal L}_{\btau}$, then
%${\btau}^{nk+1}(b_{k-1}\, l) = \btau^{nk+1}(b_{k-1})\btau^{nk+1}(l) \, \in
%{\cal L}_{\btau}$, for each $n>0$. Thus $a\,\btau^{nk+1}(l) \, \in {\cal
%L}(\btau)$ for each $n>0$, and it follows that $a\bu \,\in X$.

Conversely, suppose that $a\,\bu \, \in X_{\btau}$. We claim first that $a$ is the suffix of some
substitution word. If $a=\btau(b)$ for some $b$, we are done. If not, by Corollary \ref{almostuniquerep},  
$a\bu={\bw}\btau(\sigma \bx)$ where ${\bw}$
is the suffix of some substitution word $\btau(x_{0})$, and $\sigma
\bx$ is unique. Thus $\btau(\bu)= \bu = \sigma ({\bw})\btau(\sigma \bx)$, and by Corollary \ref{almostuniquerep} again, $\sigma (\bw) = \btau (u_{0}\ldots u_{r})$ for some $r\geq 0$, so that $a$ is the suffix of a substitution word.

Since $(X_{\btau}, \sigma)$ is minimal, $\sigma$ is onto, which means that for
each $n$ there are words $\bw_{n}$ of length $n$ such that $\bw_{n}a\bu \, \in
X_{\btau}$. Given  $N$, choose $n$ so that    $n>
|\btau^{N}(\alpha)|$ for all $\alpha$.
 Thus $\bx  = \bw_{n}a \bu =
\bw'\btau^{N}(\sigma (\by))$ where $\bw'$ is the suffix of $\btau^{N}(y_{0})$, and Corollary \ref{almostuniquerep} implies that for some $p>0$,
 $\sigma^{p+1} (\by) = \bu$ and $\bw_{n}a = \bw' \btau^{N}(\by_{[1,p]})$,
 so that $a$ is the suffix of a $\btau^{N}$ word. Taking a subsequence, we conclude that there is some $b$ such that $a$ is the suffix of $\btau^{N}(b)$ for infinitely many $N$'s. Since the sequence of suffixes of $\btau^{N}(b)$ is eventually periodic, this implies that $a$ is part of a suffix cycle.
\endproof

 %If $|{\bw}| = 1$, then ${\bw}=a$ and we are done. 
%We show that $|{\bw}|>1$ contradicts recognizability. Shifting both
%sides of Equation \ref{recog} and using part 3 of Proposition
%\ref{hosttopologyprop}, $\sigma({\bw})= \btau(\bu_{[0,r]})$ for some
%$r\geq 0$. This means that $\btau(u_{r})$ is a proper suffix of ${\bf
%w}$, and taking $W_{n}$ to be the word in $\bu$ starting at $u_{|\btau(\bu_{[0,r]})|+1}$ of length $n$, we have a contradiction to recognizability by 
%Theorem \ref{mosserecognizability}.

It follows that if $\btau$ is left proper and  suffix permutative, then $\bu$ has $|{\cal A}|$ $\sigma$-preimages, and if $\btau $ is left and right proper, then $\bu$ has one 
$\sigma$-preimage. 

Let ${\cal W}_{n} :=\{ \bw: \bw=\sigma^{k}(\btau^{n}(a)):a \in {\cal A},\,\, 0< k<|\btau^{n}(a)|\}$.
If $\bw \, \in {\cal W}_{n}$, let
$r(\bw):=|\{\alpha: \alpha\bw \, \in {\cal W}_{n}\}|$, and let
 $\sigma^{-1}\bw= \{\alpha\bw: \alpha\bw \, \in {\cal W}_{n}\}$.

\begin{proposition}
\label{sigmainjectivity}
 Let $\btau$ be reasonable, 
 with a generating fixed point $\bu$.  Suppose that $\by\neq {\bu}$, $\by\, \in
 [\bw]$ where $\bw \,\in {\cal W}_{n}$. Then $\by$ has more than
 one $\sigma$ preimage only if $r(\bw)>1$.
Conversely, suppose that there exists a sequence of words $\bw_{n}\,
\in {\cal W}_{n}$ with $|\bw_{n}|\rightarrow \infty$, with $\cap_{n}[\bw_{n}] = \by$.
If  $r({\bw}_{n})>1$ for each $n$,
 then $\by$ has more than one $\sigma$-preimage.

\end{proposition} 

{\bf Proof:} Since $\by\neq \bu$, then for some $N$, $\by\not\in \btau^{n}(X_{\btau})$ for each $n\geq N$, so that $\by= {\bw}_{n}\btau(x)$ with $\bw_{n} \in {\cal W}_{n}$. Then $\by$ has  more than one $\sigma$-preimage  only if  $r({\bf w_{n}})>1$ for each $n\geq N$.
For each $n$ there are two distinct letters
$\alpha_{n}$ and $\beta_{n}$ so that $\alpha_{n}\bw_{n}$ and
$\beta_{n}\bw_{n}$ are both in ${\cal W}_{n}$. Passing down to a
subsequence if necessary, we can assume that there are distinct
$\alpha$ and $\beta$ so that $\alpha \bw_{n}$ and $\beta \bw_{n}$ are
elements of ${\cal W}_{n}$, hence ${\cal L}_{\btau}$. Thus $\alpha \by$
and $\beta \by$ are both elements of $X_{\btau}$, so that $\by$ has more than one $\sigma$-preimage.
\endproof

Note also that  $\by =\bu$  if and only if for all large  $n$, there is a sequence $k_{n}\rightarrow \infty$. such that 
${\bf w}_{n}\,$ has the word
$\btau^{k_{n}}(u_{0})$ as a prefix.
 Proposition \ref{sigmainjectivity} suggests an algorithm for generating all non-$\sigma$ invertible points, and checking whether they are the $\btau$-fixed point. We develop this algorithm in what follows, as our aim is to classify   substitutions $\btau$ which  have a unique, non-$\sigma$-invertible point, and, when this point $\by$ is not $\btau$-fixed, we will recast  $(X_{\btau},
\sigma)$ as a substitution system for which $\by$ is fixed. 

If $\btau$ is a reasonable substitution on the alphabet ${\cal
A}=\{a_{1},\ldots, a_{n}\}$, pick ${A}_{1}\subset {\cal A}$, where
${\bf s}_{1}$, the maximal proper common suffix of $\{\btau (a)\}_{a\,
\in { A_{1}}}$, is non-empty. Thus for each $a$ in ${A_{1}}$, $\btau
(a) = {\bw}_{a}x_{a}{\bf s}_{1}$, and the set $A_{2}:=\{x_{a}:a\in
A_{1}\}$ has at least two distinct elements. Call the set $A_{2}$ the
set of ${\bf s}_{1}$-{\em predecessors} in $\{\btau(a):a\, \in
A_{1}\}$. If $\{\btau(a): a\, \in A_{2}\}$ has ${\bf s}_{2}$ as a
maximal proper common suffix, let $A_{3}$ be the set of ${\bf
s}_{2}$-predecessors in $\{\btau(a): a \, \in A_{2}\}$.  Inductively,
if $A_{k}$ is the set of ${\bf s}_{k-1}$-predecessors in $\{\btau(a):
a\,\in A_{k-1}\}$, let ${\bf s}_{k}$ be the maximal proper common
suffix of $\{\btau(a): a\, \in A_{k}\}$, and $A_{k+1}$ the set of
${\bf s}_{k}$-predecessors in $\{\btau(a): a\in A_{k}\}$.  This
process may {\em fizzle out}: $A_{k}=\emptyset $ for some $k$ - in
which case we discard this $A_{1}$ and start over with another subset of ${\cal A}$.
 The sequence $(A_{k})$ is eventually periodic. For, the
cardinality $|A_{k}|$ of the sets $A_{k}$ cannot increase. If
$|A_{k}|$ stays constant as $k$ increases, then for some $k$ and $m$,
$A_{k}=A_{k+m}$, and so ${\bf s}_{k}= {\bf s}_{k+m}$, which in turn
implies that $A_{k+1}= A_{k+m+1}$, and by induction that $A_{k}$ is
eventually periodic. If $|A_{k}|$ decreases, then repeat this
argument. Thus $({\bf s}_{k})_{k\,\in {\mathbb N}}$ is also eventually
periodic. Note that there can be more than one choice for ${
A}_{1}$ which leads to a non-trivial sequence of ${\bf s}_{k}'s$ and
$A_{k}$'s, and that there are only finitely many sequences $({\bf s}_{k})_{k\,\in {\mathbb N}}$ that can be generated, as ${A}_{1}$ ranges over all subsets of ${\cal A}$.

We illustrate this terminology with the
following examples.
\begin{enumerate}

\item 
If $\btau(a)= aab$, and $\btau(b)=abb$, if  $A_{1}={\cal A}$, then $A_{n}={\cal A}$  and  ${\bf s}_{n}=b$ for all $k\geq 2$. (In fact if $\btau$ is  any reasonable substitution on a two letter alphabet whose substitution words have a proper maximal non-empty suffix, then it is always the case that 
$({\bf s}_{k})_{k\geq 1}$ is constant.) 

\item
If $\btau(a)=acb$, $\btau(b)=aba$ and $\btau(c)=aaa$,  then the only non-fizzling choice for ${A}_{1}$ is ${ A}_{1}=\{b,c\}$. Here ${\bf s}_{1}= a$, ${ A}_{2}=\{a,b\}$, and ${\bf s}_{2} = \emptyset$ as $\btau (a)$ and $\btau (b)$ have no common suffix. Thus even this choice of ${A}_{1}$ fizzles out.

\item
If $\btau(a)=acb$, $\btau(b)=aba$ and $\btau(c)=aca$, letting ${A}_{1}=\{b,c\}$, 
then ${ A}_{k}= {A}_{1}$ and ${\bf s}_{k}={\bf s}_{2}$ for all $k\geq 2$.

\end{enumerate}

If $(Y,T)$ is a system where $\by$ is not $T$-invertible, we
call $\by$ a {\em branch point}; if $\by$ has $M$ pre-images under $T$, we
say that ${\bf y}$ is an $M$-branch point. Let us call systems $(Y,T)$
with only one branch point {\em quasi-invertible}.  If
$(X_{\btau},\sigma)$ is quasi-invertible, and the branch point is an
$M$-branch point, we'll call $\btau$ $M$-{\em quasi-invertible}.
 We use this notation in the
following lemma:

\begin{lemma}
\label{nonfixedlemma1a}
Let ${\btau}$ be a  reasonable substitution on the alphabet
${\cal A}=\{a_{1},\ldots ,a_{n}\}$, and,  given $A_{1}\subset {\cal A}$ , $({\bf s}_{m})_{m=1}^{\infty}$ and $(A_{m})_{m=1}^{\infty}$ as defined and non-fizzling. Suppose that $ {\bf s}_{k+p} = {\bf s}_{k}$ for $k\geq M$. Then 
\begin{enumerate}
\item
For 
$1\leq n\leq M+p-1$, the $\btau^{n}$-substitution words have a maximal proper common suffix 
${\bf S}_{n}=\left(\prod_{k=0}^{n-1}{\btau}^{k}({\bf s}_{n-k})\right)$, and for $n=M+kp+j$, where $k\geq 1$ and $0\leq j\leq p-1$,  
 the $\btau^{n}$-substitution words have a maximal  proper common suffix 
${\bf S}_{n}={\bf P}_{n}\,{\bf I}_{n}$, where 
 \begin{equation}
\label{yuk}{\bf P}_{n}=\left(\prod_{m=0}^{j}\btau^{m}({\bf s}_{M+j-m})\right)
\prod_{l=0}^{k-1}\left(\prod_{i=0}^{p-1}
\btau^{n-(M+(k-1-l)p+p-1-i)}({\bf s}_{M+(p-1-i)})\right)
\end{equation}
and 
\begin{equation} \label{yuk*}{\bf I}_{n}   =\left(        \prod_{m=1}^{M-1}{\btau}^{n-(M-m)}({\bf s}_{M-m})\                                                                                            \right) 
.
\end{equation}
\item
If $({\bf s}_{n})$ is eventually $k$-periodic sequence $({\bf w}_{i})_{i=1}^{k}$, then it defines $k$ branch points, with each branch point ${\by}$ satisfying an equation of the form 
\[({\bf w}_{i}\btau({\bf w}_{i+1})\btau^{2}({\bf w}_{i+2})\ldots \btau^{k-1}{\bf w}_{i-1})\btau({\by}) ={\by}\,.\]
Also any branch point ${\by}$ which is non-$\btau$-fixed arises in this fashion.
\item
The substitution $\btau$ is quasi-invertible, with branch point $\by$ non $\btau$-fixed, if and only if all $\btau$-fixed points are invertible, and there is an ${\bf s}\in {\cal A}^{+}$ such that any non-fizzling
$A_{1}\subset {\cal A}$ generates a sequence $({\bf s}_{n})$ which is
eventually the constant sequence $({\bf s})$.

\item
 The substitution $\btau$ is quasi-invertible, with branch point ${\bf u}$, a $\btau$-fixed point, if and only  ${\bf u} $ is the only $\btau$-fixed branch point, and 
% the $\btau$-substitution words 
 %have no maximal proper non-empty common suffix ${\bf s}_{1}$, and
for all non-fizzling subsets ${\cal A}_{1} \subset {\cal A}$, there is
some $N$ such that for $n>N$, there is a sequence $(k_{n})_{n\geq N}$
of natural numbers such that the maximal common proper suffix of the
$\{\btau^{n}(a):a\,\in A_{1}\}$ have $\btau^{k_{n}}(u_{0})$ as a prefix,
where $k_{n}\rightarrow \infty$.
\end{enumerate}
\end{lemma}

{\bf Proof}: 
\begin{enumerate}
\item
The proof of Equations \ref{yuk} and \ref{yuk*} follow from the definitions of $({\bf s}_{k})$, and induction.
\item
Note that if $\lim_{k\rightarrow\infty}{\bf S}_{n_k}$ exists, and $|{\bf P}_{n_k}|\rightarrow \infty$, then  
$\lim_{k\rightarrow\infty}{\bf S}_{n_k} = \lim_{k\rightarrow\infty}{\bf P}_{n_k}$, so we will assume that any sequence $({\bf s}_{k})$ is periodic.

%If $\by$ is the branch point and $\by\neq \bu$ , then $\bu$ is $\sigma$-invertible, so  by Proposition
%\ref{inverseimagesofu}, working with a higher power of $\btau$ if necessary, we can assume that the the $\btau$-substitution words have a maximal
%common suffix ${\bf s}_{1}$, and ${\bf s}_{1}$ has to be proper, or else
%$(X_{\btau},\sigma)$ is periodic and so invertible.
In this case  
Equation \ref{yuk} reduces to
\begin{equation}
\label{lesseryuk}
{\bf S}_{np}={\bf P}_{np}=\left(\prod_{m=0}^{p-1}\btau^{m}{\bf s}_{p-m}\right)
\prod_{l=0}^{n-1}\left(\prod_{i=0}^{p-1}
\btau^{lp+i}({\bf s}_{p-i})\right)
\end{equation}

Here ${\bf S}_{(n+1)p}$ is a prefix of ${\bf S}_{np}$ for $n\geq 1$, 
and 
%if the $\btau$-substitution words start with $u_{0}$, then $u_{0} \,\in {\bf
%S}_{np}$ for $n>1$, so that
if  $|{\bf S}_{np}|\rightarrow \infty$,
then
$\by:=\lim_{n\rightarrow \infty}{\bf S}_{np}$ exists,  and  by Proposition \ref{sigmainjectivity}  has at  least 2 $\sigma$-preimages. Note that $\by$ satisfies the equation 
\begin{equation}\label{theone}\left(\prod_{m=0}^{p-1}\btau^{m}{\bf s}_{p-m}\right)\,\btau^{p}({\bf y})= {\bf y}\,.\end{equation}
%If the sequence $({\bf s}_{n})$ is not eventually fixed -ie if $p>1$, a formula similar to Equation \ref{lesseryuk}
%for $({\bf S}_{np+k})_{n\geq 0}$ can be obtained for  $0<k<p$. These words start with 
%$\prod_{m=0}^{k-1}\btau^{m}({\bf s}_{1+k-m}) $, leading to another branch point, satisfying a different equation to Equation \ref{theone}. This contradicts quasi-invertibility.
A similar argument applied to the sequence $({\bf S}_{(n+1)p+i})_{n}$ yields the other branch points.

 Conversely, if $\by$ is a non-$\btau$-fixed $k$-branch point, then
 there are $k$ letters $\{a_{1}^{n},a_{2}^{n},\ldots a_{k}^{n}\}$
 whose $\btau^{n}$ words share a maximal common proper suffix ${\bf
 T}_{n}$ and such that $\by\in [{\bf T}_{n}]$. Moving to a subsequence
 if necessary, we can assume that these $k$ letters are independent of
 $n$ and so if $A_{1}=\{a_{1}, \ldots a_{k}\}$, the sequence $({\bf
 T}_{n})$ is identical to the sequence $({\bf S}_{n})$ generated by
 $A_{1}$. The result follows.

% the $\btau$-substitution words (we can assume $L=1$)
%$\btau(a_{i})$  have a maximal proper non-empty common suffix ${\bf s}_{1}$, and $p=1$,  then ${\bf u}$ is not a branch point. Note that by Proposition \ref{sigmainjectivity} a branch point ${\bf y}$ must be ${\bf y}=\lim_{n\rightarrow \infty}{\bf S}_{n}$ for the sequence $({\bf S}_{n})$ generated as in Equation \ref{lesseryuk} by some appropriate $A_{1}$. As in Equation \ref{theone}, the point ${\bf y}$ thus has to satisfy ${\bf w}\btau^{k}(y)={\bf y}$ for some word ${\bf w}$ and some natural $k$, both generated by the sequence $({\bf s}_{k})$. If all sequences $({\bf s}_{k})$ generated by any choice of $A_{1}$ are ultimately identical and constant, then only one branch point ${\bf y}$ can be generated.

\item

Statements 3  and 4  now follow from Part 2.
\endproof

\end{enumerate}

\begin{lemma}
\label{nonfixedlemma2a}
Suppose that the reasonable $\btau$ is defined as  $\btau (a_{i})= {\bf p}_{i} x_{i} {\bf s}_{1}$, where
$\{x_{1},\ldots, x_{n}\}$ has at least two distinct elements, and where if ${ A}_{1}= {\cal A}$, then  
${\bf s}_{k}= {\bf s}_{1}$ for $k\geq 1$.  Define
the substitution ${\btau}^{*}$ as ${\btau}^{*}(a_{i}) = {\bf
s}_{1}{\bf p}_{i}x_{i}$.
\begin{enumerate}
\item
For each $n\geq 0$, and each $\bw \, \in {\cal L}_{\btau}$,
${\btau}^{*}({\btau}^{n}(\bw)){\bf s}_{1} = {\bf s}_{1} \btau^{n+1}(\bw)\,.$
\item
 If ${\bf S}_{n}$ is as in Lemma \ref{nonfixedlemma1a}, then for each $n \geq 2$,
${\bf S}_{n}= \btau^{*}({\bf S}_{n-1})\,{\bf s}_{1}\, .$
\item
${\bf S}_{n}\subset ({\btau}^{*})^{n-1}({\bf s}_{1})$
for each $n> 1$.
\end{enumerate}

\end{lemma}
{\bf Proof}:
\begin{enumerate}
\item
 %We prove this statement only for $\bw =a_{1}$; the proof for other one-letter words is the same, and the proof for longer words follows by concatenation. 
Note that 
$\btau^{*}(a)\, {\bf s}_{1}= {\bf s}_{1}{\bf p}_{a}x_{a}{\bf s}_{1} = {\bf s}_{1}{\btau}(a)\, $ for any $a\in {\cal A}$.  The general proof for longer words follows by concatenation.
% If the assertion is true for $n$, then 
%\[\btau^{*}(\btau^{n+1}(a_{1})){\bf s}_{1} = \btau^{*}(\btau^{n}(\btau (a_{1}))) {\bf s}_{1}
%= \btau^{*}(\btau^{n}({\bf p}_{1} x_{1} {\bf s}_{1}  )) {\bf s}_{1} = {\bf s}_{1} \btau^{n+1}(\btau(a_{1})) = {\bf s}_{1} \btau^{n+2}(a_{1})\, .\]
\item
 First by Lemma \ref{nonfixedlemma1a} and Part 1 of this lemma, 
${\bf S}_{2}= {\bf s}_{1}\btau ({\bf s}_{1})= \btau^{*}({\bf s}_{1}){\bf s}_{1}$. Assuming that ${\bf S}_{n}= \btau^{*}({\bf S}_{n-1}){\bf s}_{1}$, we have 
\[{\bf S}_{n+1} = {\bf S}_{n} \btau^{n}({\bf s}_{1})
\stackrel{IH}{=} \btau^{*}({\bf S}_{n-1}){\bf s}_{1} \btau^{n}({\bf s}_{1})
\stackrel{Part 1}{=} 
\btau^{*}({\bf S}_{n-1}){\btau}^{*}(\btau^{n-1}({\bf s}_{1})){\bf s}_{1}
= \btau^{*}({\bf S}_{n})  {\bf s}_{1}.\]
\item
The case $n=2$ is clear.
Assuming ${\bf S}_{n}\subset ({\btau}^{*})^{n-1}({\bf s}_{1}),$
\[{\bf S}_{n+1} \stackrel{Part 2}{=} {\btau}^{*}({\bf S}_{n}) {\bf s}_{1}
\stackrel{IH}{\subset} {\btau}^{*}(({\btau}^{*})^{n-1})({\bf s}_{1})) {\bf s}_{1}
= ({\btau}^{*})^{n}({\bf s}_{1}) \, {\bf s}_{1}\subset ({\btau}^{*})^{n}({\bf s}_{1})\,.\]
\endproof\end{enumerate}

%Lemma \ref{nonfixedlemma2a}    can be modified if $M>1$:  instead of considering $\btau$, work with $\btau^{M}$, and as in the lemma, find $\btau^{*}$ such that $(X_{\btau^{*}},\sigma)=(X_{\btau^{N}},\sigma)$; now use the fact that $(X_{\btau^{N}},\sigma)=(X_{\btau}, \sigma)$.

\begin{corollary}
\label{nonfixedcorollary}

If the left proper, reasonable substitution subshift 
$(X_{\btau},\sigma)$ is quasi-invertible,    
where $\tau$ is defined on ${\cal A}$ with 
branch point $\by\neq \bu$, then some power of $\btau$ is right proper, and there exists a left proper, recognizable, minimal, aperiodic, injective, quasi-invertible 
substitution subshift with  
$\btau^{*}(\by) = \by$, and $(X_{\btau},\sigma)=(X_{\btau^{*}}, \sigma)$.
\end{corollary}

{\bf Proof}: Since the fixed point $\bu$ is $\sigma$-invertible, Proposition  \ref{inverseimagesofu} tells us that some $\btau$ (or some power of $\btau$) is right proper.
If $\by$ is the  branch point for $\btau$,
then Lemma \ref{nonfixedlemma1a} tells us that there exists an ${\bf s}$ such that if $A$ is any subset of  ${\cal A}$, it  generates an eventually fixed sequence $({\bf s}_{k})$ where  ${\bf s}_{k}= {\bf s}$ for $k\geq n_{A}$ (No $A$ fizzles out as $\btau$ is left proper). 
 Since there are finitely many such sequences we can assume, taking a power if necessary,  that $n_{A}=1$.  Now Lemma \ref{nonfixedlemma2a} can be applied, taking 
$\btau^{*}$ as defined. The branch point is  $\by:=\lim_{n\rightarrow \infty}{\bf S}_{n}$; that it
exists and has at least two $\sigma$-preimages follows from Lemma
\ref{nonfixedlemma1a}. That it is ${\btau^{*}}$-fixed follows from
Part 3 of Lemma \ref{nonfixedlemma2a}. Thus $X_{\btau^{*}}\subset
X_{\btau}$, and by minimality, this inclusion is an equality.
\endproof

{\bf Examples}
All substitutions are reasonable.
\begin{enumerate} \setcounter{enumi}{3}

\item(Example 1)
If $\btau(a)= aab$, and $\btau(b)=abb$, then  the unique branch point $\by$ satisfies $b\btau (\by)= \by$.

\item(Example 2)
If $\btau(a)=acb$, $\btau(b)=aba$ and $\btau(c)=aaa$,  then only $\bu$ is a 2-branch point.

\item(Example 3)
If $\btau(a)=acb$, $\btau(b)=aba$ and $\btau(c)=aca$, then $\{a,b\}$ is a suffix cycle, so $\bu$ is a 2-branch point. Letting ${A}_{1}
=\{b,c\}$, we obtain the (only other) 2-branch point $\by$ satisfying $a\btau(\by)=\by$, and $\by\neq \bu$.

\item If $\btau(a)=abc$, $\btau(b)=aacc$, and $\btau(c)=abcc$, then $\bu$ is $\sigma$-invertible, and  if 
$A_{1}=\{b,c\}$,  ${\bf s}_{2k+1}=cc$ and ${\bf s}_{2k}=c$ for  $k\geq 1$, and no other choice of $A_{1}$ generates any other eventually different sequence $({\bf s}_{k})$.
Thus $(X_{\btau},\sigma)$ has two branch points, $\by$ and $\by^{*}$, satisfying $\by=c\btau(cc)\btau(\by)$ and $\by^{*}= 
cc\btau(c)\btau(\by^{*})$, and these two points are distinct.

%\item
%If $\btau(a)=abd$, $\btau(b)=abcd$, $\btau(c)=accd$ and $\btau(d)=aacd$, then for $k\geq 2$, ${\bf e}_{k}{\bf s}_{1}=cd$,  $A_{k}=\{b,c\}$,   so that $(X_{\btau},\sigma)$ has exactly  one branch point, satisfying $cd\btau(\by)=\by$.

%\item There are non suffix-permutative substitutions where $\bu$ is the only non invertible element. The substitution defined on ${\cal A}=\{a,\, b, \, c, \, d\}$ defined by 
%$\btau(a)= ad$, $ \btau (b) = a$, $\btau (c) = d$ and $\btau (d) = bc$ is recognizable (see \cite{Mo1}), injective, and minimal. It is not suffix permutative, yet only has one point with multiple $\sigma$-preimages: only $c$ and $d$ are part of a suffix cycle, so $\bu$ has 2 preimages. Note that $\btau^{2}(b) = \btau^{2}(d)$, and $\btau (c)$ is a proper suffix of $\btau (a)$. Thus there are no other elements with multiple $\sigma$-preimages.

\item 
If  $\btau(a)=bbad,$ $\btau(b)=ab$, $\btau(c)=ad$, and $\btau(d) = dac$, then $\btau$ has 3 fixed points, with $\btau^{\infty}(a)$ a 3-branch point, $\btau^{\infty}(b)$ a 2-branch point, and where $\btau^{\infty}(d)$ is $\sigma$-invertible.

%is  recognizable: although $\btau(c)$ is a proper suffix of $\btau(a)$, $a$ is always followed by $b$ or $d$ in $X$, while $c$ is only followed by $a$. So $\btau(a)$ is always followed by $\btau(b)$ or $\btau(d)$, and $\btau(c)$ is always followed by $\btau(a)$. Hence long words satisfying part 2 of Theorem \ref{mosserecognizability} do not exist.
%As $\{b,c,d\}$ belong to suffix cycles, $\bu$ has three $\sigma$-preimages.
If $A_{1}= \{a,b\}$ then the words $\{\btau^{n}(a), \btau^{n}(b)\}$ have $\btau^{n}(b)$ as a maximal common suffix, and these words also converge to $\bu$.

%Now $\btau^{n}(b)$ always has $b$ as a suffix, and no other substitution words do. Also,$\btau^{2n+1}(a)$
%always ends with $d$, $\btau^{2n}(a)$ always ends with $c$, while the
%opposite is true for the $\btau^{n}(d)$ substitution words of
%$a$. $\btau(c)$ is a proper suffix of $\btau(a)$.  Thus the only time
%two $\btau^{n}$-words have a common suffix, one of the $\btau^{n}$-words is a suffix of another, 

%\item

%If $\btau(a)=aad$, $\btau(b)=abd$, $\btau(c)=
%ac$, and $\btau(d)=abc$. All ${\btau}^{2}$-substitution words have  $c$ as a suffix,  so $\bu$ is $\sigma$-invertible.
%Also, ${\bf e}_{k}=dab$ for $k\geq 2$, and $A_{k}=\{a,b\}$ for all $k$. Working with ${\btau}^{4}$, there is  only  one  branch point $\by$; it satisfies 
%$dabc \, \btau(\by)=\by$.

%\item If $\btau(a)=abd$, $\btau(b)=aad$, $\btau(c)=
%accd$, and $\btau(d)=adcd$. It has one branch point $\by=d\btau(\by)$; it also has another satisfying $cd \btau(\by)=\by$, arising from ${\cal A}^{*}=\{c,d\}$.

\end{enumerate}

Unilateral recognizability is less amenable than  bilateral
recognizability, which all aperiodic primitive substitutions possess,
and for whom the sequence of partitions $({\cal P}_{n})$ defined by
\[{\cal P}_{n}^{'}
:=
\{[\sigma^{k}(\btau^{n}(a))]: a \, \in {\cal A},\,\,\, 0 \leq k < |\tau^{n}(a)|\}
\]
is a nested sequence which spans the topology of $(X^{\mathbb Z}_{\btau},\sigma)$ (see  Proposition 14, \cite{dhs}). In the one sided subshift $(X_{\btau}, \sigma)$, the sets in 
${\cal P}_{n}^{'}$ need not be disjoint. For example, if $\btau$ has a branch point $\by\neq \bu$, then, for each $n$, there are letters $a$, $b$, and numbers $n_{1}$ and $n_{2}$ such that $|\btau^{n}(a)|-n_{1} = |\btau^{n}(b)|-n_{2}$ and  
$\sigma^{n_{1}}(\btau^{n}(a))=\sigma^{n_{2}}(\btau^{n}(b))$, while 
$\sigma^{n_{1}-1}(\btau^{n}(a))\neq\sigma^{n_{2}-1}(\btau^{n}(b))$. If $\btau$ is such that $\bu$ is the only branch point in $X_{\btau}$, or even if the branch point lies in $\btau(X_{\btau})$,  this does not happen.

% To get around
%this  we   define
%\[
%{\cal W}_{n}:= \{ \bw: \bw = \sigma^{k}(\btau^{n}(a)): a \, \in {\cal A} \mbox{ and } 0<k< |\btau^{n}(a)|
 %\}\, .\]
%The results in Corollary \ref{almostuniquerep} can be extended to tell
%us that for each $n\geq 1$, every element $\by\,\in X\backslash \btau^{n}
%(X)$ has a unique representation $\by=\bw_{n}\bx_{n}$ with $\bw_{n} \,
%\in {\cal W}_{n}$ and $\bx_{n}\, \in \btau^{n}(X)$. Given $\bw \, \in {\cal W}_{n}$, define 
%\[ [\bw]^{*}:=  [\bw]\cap \sigma^{-|\bw|}\,\btau^{n}(X) =  \{ \by \, \in X: 
%\by = \bw \bx \mbox{ with } \bx \, \in {\btau}^{n}(X) \}. \]
% We assume now that $\btau$ is primitive, so that $\lim_{n\rightarrow \infty}|\btau^{n}(a)|= \infty$ for all $a\, \in {\cal A}$. Define 
%\[
%{\cal P}_{n}:= \{\btau^{n}(a): a\, \in {\cal A} \} \cup  \{ [\bw]^{*}: \bw \, \in {\cal W}_{n}\}. \]
%Call the set $\{\btau^{n}(a): a\, \in {\cal A} \} $ the {\em base} of ${\cal P}_{n}$. Sets $[\bw]^{*}$ will either be a single `level', $[\bw]^{*}= [\sigma^{k}(\btau^{n}(a))]$ in $X_{\btau}$, or a union 
%\[[\bw]^{*} = \cup_{i=1}^{k} [\sigma^{n_{i}}(\btau^{n}(a_{i}))],\]
%where $|\btau^{n}(a_{i})|-n_{i} = |\btau^{n}(a_{j})|-n_{j}$ for each $i$, $j\, \in \{1,2,\ldots k\}$. This happens when the words $\{\btau^{n}(a_{j})\}_{j=1}^{n}$ have a common suffix of length 
%$|\btau^{n}(a_{1})|-n_{1}$.

\begin{corollary}
\label{clopenpartition}
Let $\btau$ be  left proper and reasonable, and suppose that $\bu$ is the only branch point in $X_{\btau}$. 

Then
\begin{enumerate}
\item
For every $n>0$, ${\cal P}_{n}$ is a clopen partition, and the sequence of bases is decreasing. 
\item For each $n$, ${\cal P}_{n+1}$ is a refinement of ${\cal P}_{n}$: i.e.
every element of ${\cal P}_{n+1}$ is contained in an element of ${\cal
P}_{n}$.
\item
The intersection of the bases of ${\cal P}_{n}$ consists of a unique point, and the sequence $({\cal P }_{n} )$ spans the topology of $(X_{\btau},\sigma)$.

\end{enumerate}
\end{corollary}

{\bf Proof:}
\begin{enumerate}
\item
We show that ${\cal P}_{1}$ is clopen, First suppose that $({\bf
y}_{n})$ is a sequence in $\btau [a]$ and ${\bf y}_{n}\rightarrow
{\bf y}$. Since $\btau(X_{\btau})$ is closed, by Proposition
\ref{hosttopologyprop}, then ${\bf y}\, \in \btau(X_{\btau})$ If ${\bf
y}_{n}=\btau({\bf x}_{n})$ and ${\bf x}$ is a limit point of $({\bf
x}_{n})$, then ${\bf x}\,\in [a]$ and ${\bf y}=\btau({\bf x})\subset
[\btau(a)].$ Using Corollary \ref{almostuniquerep}, the sets
$\{\btau[a]: a\, \in {\cal A}\}$ are disjoint, and their union,
$\btau(X_{\btau})$, is a clopen set, so they are also.  If $0<k<|\btau (a)|$
and $({\bf y}_{n})$ is a sequence of elements in
$[\sigma^{k}(\btau(a))]$, then ${\bf y}_{n}=\sigma^{k}({\bf x}_{n})$
where ${\bf x}_{n}\,\in \btau([a])$. If ${\bf x}$ is a limit point of
$({\bf x}_{n})$ then it must also be in $\btau[a]$, and ${\bf y} =
\sigma^{k}({\bf x})\in \sigma^{k}(\btau[a])$. So the sets
$\sigma^{k}(\btau[a])$ are closed, and if the branch point is the
fixed point $u$, then these sets are also disjoint, and so clopen.

%To show that they are disjoint
%suppose that $\by\,\in [\bw]_{*}\cap [{\bw'}]_{*}$ where $\bw =
%\sigma^{k}(\tau^{n}(a))$ and ${\bw'}=\sigma^{l}(\btau^{n}(b))$. If
%$|\bw| <|\bw'|$ (or $|\bw| >|\bw'|$), then $\sigma^{|\bw|}(\by)$
%( or $\sigma^{|\bw'|}(\by)$) has two representations of the form
%$\sigma^{j}(\btau(\bx))$, a contradiction, using Corollary
%\ref{almostuniquerep}. So $|\bw|=|\bw'|$, which means that $\bw$ is a
%common suffix of $\btau^{n}(a)$ and $\btau^{n}(b)$. Since
%$\btau^{n}(a)\neq \btau^{n}(b)$, these two words have a maximal proper
%common suffix, say $\sigma^{k'}(\tau^{n}(a))=\sigma^{l'}(\btau^{n}(b))$, and here one of $k'$ or $l'$ must be positive - say $k'$. Thus $\sigma^{-(k-k')}(\by)=\sigma^{-(l-l')}(\by)$ has two preimages, so by assumption 
%$\bu= \sigma^{-(k-k')}\by\,\in \sigma^{-(k-k')}[\sigma^{k}(\btau^{n}(a))]_{*}
%=[\sigma^{k'}\btau^{n}(a)]_{*}, $ once again a contradiction to Corollary \ref{almostuniquerep}.
\item
 Take an element ${\bf p} \in {\cal P}_{n+1}$: we need to
show that $[\bp]\subset [\bp']$ where $[\bp'] \in {\cal P}_{n}$.
Suppose that 
$\bp =
\sigma^{k}(\btau^{n+1}(a))$, 
where $0<k< |\btau^{n+1}(a)|$. 
As in the proof of Proposition 14, in \cite{dhs}, write 
$\btau^{n+1}(a)= \btau^{n}(a_{1})\,\btau^{n}(a_{2})\, \ldots
\btau^{n}(a_{m})$ where $\btau(a)= a_{1}\, a_{2} \ldots a_{m}$. For
some $j$, $|\btau^{n}(a_{1}\ldots a_{j})|\leq k <
|\btau^{n}(a_{1}\ldots a_{j+1})|$, so let $l:=k -
|\btau^{n}(a_{1}\ldots a_{j})|$. Then
$\sigma^{k}(\btau^{n+1}([a]))\subset \sigma^{l}\btau^{n}[a']$, where 
$a'=a_{j+1}$. Let $\bp':= \sigma^{l}\btau^{n}[a']$.

 \item
 To see that $({\cal
P}_{n})$ span the topology of $(X_{\btau},\sigma)$, we imitate the proof of
Proposition 14, \cite{dhs}. Given $m$ positive, we show that for all  $n$
large,  each element ${\bf p}$ of
${\cal P}_{n}$ is contained in some ${\bf c}$ where ${\bf c}$ is a word of 
length $m+1$. Suppose that all $\btau$-substitution words start with
$l$. Write $l_{n} = |\btau^{n-1}(l)|$. Choose $n$ so that $l_{n}>m$
(Minimality implies that $\lim_{n\rightarrow \infty}l_{n} = \infty$).

Fix $\bp \,\in {\cal P}_{n}$, and suppose that $\bp = \sigma^{k}(\btau^{n}(a))$ for some  $a\,\in {\cal A}$ and $k\,\in [0, |\btau^{n}(a)|)$. If $\by \,\in
[\bp]$ then there is some $\bx\,\in [a]$ with
$\by=\sigma^{k}\btau^{n}(\bx)$. Since $\btau$ is left proper,
${\btau}^{n}(a)\btau^{n-1}(l)$ is a prefix of ${\btau}^{n}(\bx)$, and
so $\by \in [\sigma^{k}{\btau}^{n}(a)\btau^{n-1}(l)], $ and this last
cylinder set has length greater than $m$. So, $\by$ starts with a
block of length $m$ depending only on $\bp$ and  not on $\by$.  \endproof
\end{enumerate}

%Reem: The above is just a simple generalization of DHS. So maybe delete it later.

\section{Bratteli Diagrams}\label{Bratteli}
A {\em Bratteli diagram} ${\cal B}=({\cal V},{\cal E})$ is an infinite
directed graph with {\em vertex set} ${\cal V} = \bigsqcup_{n=
0}^{\infty} {\cal V}_{n}$ and {\em edge set} ${\cal E} =
\bigsqcup_{n=1}^{\infty}{\cal E}_{n}$, where all ${\cal V}_{n}$'s and
${\cal E}_{n}$'s are finite, ${\cal V}_{0}= \{v_{0}\}$, and, if $x$ is
an edge in ${\cal E}_{n}$, the {\em source} $s(x)$ of $x$ lies in
${\cal V}_{n-1}$ and the {\em range} $r(x)$ of $x$ lies in ${\cal
V}_{n}$. We assume that $s^{-1}(v)\neq \emptyset$ for each $v\in {\cal
V}$ and $r^{-1}(v) \neq \emptyset $ for all $v \, \in {\cal V}
\backslash v_{0}$.  We will use $x, y\ldots$ when referring to edges,
and $a,b, \ldots$ when referring to vertices.

A finite set of edges $\{x_{n+k}\}_{k=1}^{K}$, with
$s(x_{n+k+1})=r(x_{n+k})$ for $1\leq k\leq K-1$, is called a {\em path
} from $s(x_{n+1})$ to $r(x_{n+K})$. Similarly an {\em infinite path}
in ${\cal B}$ is a sequence ${\bx}=(x_{n})_{n=\infty}^{1}$, with
$x_{n}\,\in {\cal E}_{n}$ for $n\geq 1$, and $s(x_{n+1})=r(x_{n})$ for
$n\geq 1$.  We  write
$
\ldots a_{n+1}\stackrel{x_{n+1}}{\leftarrow}a_{n}
\stackrel{x_{n}}{\leftarrow}\ldots \stackrel{x_{2}}{\leftarrow} a_{1}
\stackrel{x_{1}}{\leftarrow} v_{0}
$
where $r(x_{n})=a_{n}$, or $({\bf a},{\bx}) =
((a_{n},x_{n}))_{n\geq 1}$, when referring to
an element ${\bx}$ in $X_{\cal B}$.  The set of all infinite paths in ${\cal B}$
will be denoted $X_{\cal B}$ (a subset of $ \Pi_{n\geq 1} {\cal E}_{n}
$), and  $X_{\cal B}$ is endowed with the topology induced from the product
topology on $\Pi_{n\geq 1} {\cal E}_{n}$. Thus $X_{\cal B}$ is a
compact metric space.

Two Bratteli diagrams ${\cal B}=({\cal V}, {\cal E})$ and ${\cal
B'}=({\cal V'}, {\cal E'})$ are {\em isomorphic} if there exists a
pair of bijections $f_{{\cal V}}: {\cal V}\rightarrow {\cal V'}$ and
$f_{\cal E}:{\cal E}\rightarrow {\cal E'}$ satisfying $f_{\cal V}(a)
\, \in {\cal V}_{n}'$ if $a \, \in {\cal V}_{n}$, and $s(f_{\cal
E}(x)) = f_{\cal V}(s(x))$, $r(f_{\cal E}(x)) = f_{\cal V}(r(x))$
whenever $x\, \in {\cal E}$.  Let $(n_{k})_{k=0}^{\infty}$ be a
sequence of increasing integers with $n_{0}=0$. Then ${\cal B}'
=({\cal V}',{\cal E}')$ is a {\em telescoping} of ${\cal B}= ({\cal
V},{\cal E})$ if ${\cal V}_{k}'=V_{n_k}$ (with the vertex $v \, \in
{\cal V}_{n_k}$ labelled as $v' \, \in {\cal V}_{k}'$), and the number
of edges from $v_{k}'\, \in {\cal V}_{k}'$ to $v_{k+1}'\, \in {\cal
V}_{k+1}'$ is the number of paths from $v_{k}\, \in {\cal V}_{n_{k}}$
to $v_{k+1}\, \in {\cal V}_{n_{k+1}}$.  Conversely, we can perform a
{\em splitting} by introducing a new level between two consecutive
levels ${\cal V}_{n-1}$ and ${\cal V}_{n}$, so that the number of new
vertices equals the edges in ${\cal E}_{n}$, and each vertex in the
new level is the source and range of exactly one edge.  We consider
two Bratteli diagrams ${\cal B}$ and ${\cal B}'$ {\em equivalent} if
${\cal B'}$ can be obtained from ${\cal B}$ by isomorphism,
telescoping and splitting. Thus when we talk about a Bratteli diagram
we are talking about an equivalence class of diagrams.

We say that ${\cal B}$ is {\em simple} if there
exists a telescoping ${\cal B'}=({\cal V'},{\cal E'})$ of ${\cal B}$
so that, for any $a \, \in {\cal V}'_{n}$ and $b\, \in {\cal
V}'_{n+1}$, there is at least one edge from $a$ to $b$. If ${\bx}=
(x_{n})_{n=\infty}^{1}$ and ${\bx'}=(x_{n}')_{n=\infty}^{1}$ are
two elements in $X_{\cal B}$, we write ${\bx} \sim {\bx'}$ if the
tails of ${\bx}$ and ${\bx'}$ are equal. It follows that $\sim$ is
an equivalence relation, and if ${\cal B}$ is simple, then each
equivalence class for $\sim$ is dense in $X_{\cal B}$, and $X_{\cal
B}$ has no isolated points, making it a Cantor space.

\subsubsection{Ordering  $X_{\cal B}$}
Let $n\geq 1$. For each $a \, \in {\cal V}_{n}$, let ${\cal E}_{n}(a)=
\{x \, \in {\cal E}_{n}: r(x)=a\}.$ Say ${\cal B}$ is {\em ordered }
if there is a linear order $\geq$ on each ${\cal E}_{n} (a)$; elements
of ${\cal E}_n (a)$ will then be labelled $1,2, \ldots$ according to
their order. If $a\, \in {\cal V}\backslash \{v_{0}\}$, define $|a|:=
|{\cal E}_{n}(a)|$, so that ${\cal E}_{n}(a) = \{1,2, \ldots |a|\}$.

The linear order on
edges in each ${\cal E}_{n}(a)$ induces a partial ordering on paths from
${\cal V}_{m}$ to ${\cal V}_{n}$: the two paths ${\bx}=
a_{n}\stackrel{x_{n}}{\leftarrow}a_{n-1}
\stackrel{x_{n-1}}{\leftarrow}\ldots
\stackrel{x_{m+1}}{\leftarrow}a_{m}$ and ${\bf
x'}=a_{n}'\stackrel{x_{n}'}{\leftarrow}a_{n-1}'
\stackrel{x_{n-1}'}{\leftarrow}\ldots
\stackrel{x_{m+1}'}{\leftarrow}a_{m}'$ from ${\cal E}_{m}$ to ${\cal
E}_{n}$ are comparable with ${\bx}<{\bx'}$ 
if there is some $k\, \in [m+1,n]$ with $x_{k}<x_{k}'$ and
$x_{j}=x_{j}'$ for $k+1\leq j\leq n$.

Finally, two elements ${\bx},\, {\bx'}\,\in X_{\cal B}$ are {\em comparable} with
${\bx}<{\bx'}$ if there is a $k$ such that $x_{n}=x_{n}'$
for all $n>k$,  and $ x_{k}<x_{k}'$. Thus each
equivalence class for $\sim$ is ordered.
There is the obvious notion of
ordered isomorphism of two ordered Bratteli diagrams ${\cal B},$
${\cal B'}$: the isomorphism between ${\cal B}$ and ${\cal B'}$ also
has to satisfy $f_{\cal E}(x) \leq f_{\cal E}(y)$ if $x\leq y$. If
${\cal B}'$ is a telescoping of the ordered Bratteli diagram ${\cal
B}$, then the order induced on ${\cal B}'$ from the order on ${\cal B}$
makes ${\cal B}'$ an ordered Bratteli diagram.  We say that the ordered
Bratteli diagrams ${\cal B},$ ${\cal B'}$ are {\em equivalent} if
${\cal B'}$ is the image of ${\cal B}$ by telescoping and order
isomorphism.

An infinite path is {\em maximal} ({\em minimal}) if all the edges
making up the path are maximal (minimal).  If ${\bx}=
(x_{n})_{n=\infty}^{1}\, $ is not maximal, let $k$ be the smallest
integer such that $x_{k}\, \in {\cal E}_{k}(a_{k})$ is not a maximal
edge, and let $y_{k}$ be the successor of $x_{k}$ in ${\cal
E}_{k}(a_{k})$. Then the {\em successor} ${\cal V}_{\cal B} ({\bx}) $
is defined to be ${\cal V}_{\cal B}({\bx})= \ldots x_{k+2}\,x_{k+1}\,
y_{k}\, 1 \ldots 1$, where $1\ldots 1$ is the unique minimal path starting at $v_{0}$ and ending at the source of $y_{k}$. Similarly, every non-minimal path has a unique
{\em predecessor}.  Let $X_{\min}$ $(X_{\max}) \subset X_{\cal B}$ be
defined as the set of minimal (maximal) elements of $X_{\cal B}$. By
compactness, these sets are non empty.  Simple ordered Bratteli
diagrams which have a unique minimal and maximal element (called
${\bx}_{\min}$ and ${\bx}_{\max}$ respectively) are called {\em
proper}, and those with a unique minimal element are called {\em
semi-proper}.  If ${\cal B} $ is semi-proper, then $V_{\cal B}$ can be
extended to a continuous surjection on $X_{\cal B}$ by setting
$V_{\cal B} ({\bx}_{\max})= {\bx}_{\min}$, which is a homeomorphism if
${\cal B}$ is proper.  We call $(X_{\cal B},V_{\cal B})$ a {\em
Bratteli-Vershik} or {\em adic} system. Note that $(X_{\cal B},V_{\cal
B})$ is a minimal system, since $V_{\cal B}$ orbits are equivalence
classes for $\sim$. Let us say that two Cantor systems
$(X_{i},T_{i},{\bx}_{i})$, $i=1,2$ are {\em pointedly isomorphic} if
there exists a homeomorphism $f:X_{1}\rightarrow X_{2}$ with $f\circ
T_{1}= T_{2}\circ f$ and $f(x_{1})=x_{2}$.  The Bratteli- Vershik
system associated to an equivalence class of ordered semi-proper
Bratteli diagrams is well defined up to pointed isomorphism. This was
proved in Section 4, \cite{hps} for proper Bratteli diagrams; the
proof is similar for semi-proper Bratteli diagrams:

\begin{proposition}
\label{putnam}
Let ${\cal B}$ and ${\cal B}'$ be semi-proper ordered Bratteli
diagrams, with the same number of maximal elements. Then ${\cal B}$ is equivalent to ${\cal B}'$ if and only if
$(X_{\cal B}, V_{\cal B}, {\bx_{\min}})$ is pointedly isomorphic to
$(X_{\cal B'}, V_{\cal B'}, {\bx_{\min}^{'}})$.
\hfill\endproof
\end{proposition}

\begin{figure}[h]
\centerline{\includegraphics[scale=1.0]{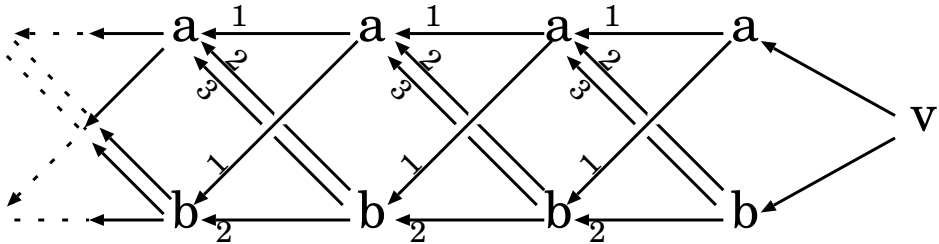}}
\caption{{\footnotesize \label{Brattelidiagram}The Bratteli diagram associated with the substitution
$\btau(a) =abb,\,\,\, \btau(b)=ab\,$.}}
\end{figure}

%\subsubsection{The Bratteli diagram associated with a proper substitution}
Let $\btau$ be a primitive substitution on ${\cal A}$.  The {\em Bratteli
diagram associated with $\btau$} has vertex sets ${\cal V}_{n} ={\cal
A},$ for each $n\geq 1$.  There is exactly one edge from each vertex
in ${\cal V}_{1}$ to $v_{0}$.  If $\btau(a) = a_{1}a_{2}
\ldots a_{n}$,
then there is an edge from vertex $b$ in ${\cal V}_{n-1}$ to vertex $a \,\in {\cal V}_{n}$, and it is labelled $i$, if and only if $a_{i}=b$.
If $1 \leq m<n$, then the number of paths in ${\cal B}$
from $(a,n) $ to $(b,m)$ is the number of occurrences of $b$ in
${\btau}^{  (n-m)}(a).$ Primitivity implies that there is a
positive $k$ such that for any two letters $a$ and $b$, there is at
least one path from $(a,n+k)$ to $(b,n)$. The Bratteli diagram ${\cal
B}' =({\cal V}',{\cal E}')$ for $\btau^{  k}$ is the telescoping
of the Bratteli diagram ${\cal B}= ({\cal V},{\cal E})$ for $\btau$,
with ${\cal V}_{n}'= {\cal V}_{nk}$ for $n>1$ and ${\cal V}_{1}' =
{\cal V}_{1}$. An example of a Bratteli diagram  is illustrated in Figure \ref{Brattelidiagram}.

The connection between Bratteli-Vershik systems and one sided left proper reasonable
substitutions is given by the next result. It is the appropriate
generalization of \cite[Prop 20]{dhs} to one sided substitution
systems.

\begin{theorem}
\label{quasiinvertibleleftproperareadic}

If $\btau$ is left proper, reasonable and  quasi-invertible, with fixed point $\bu$, there
exists a semi-proper, stationary Bratteli diagram ${\cal B}$ which is semi-proper, and
stationary, such that $(X_{\cal B}, V_{\cal B})$ is topologically
conjugate to $(X_{\btau}, \sigma)$.

\end{theorem}

{\bf Proof}: If $\bu$ is the branch point, let ${\cal B}$ be the
Bratteli diagram associated with $\btau$. Since $\btau$ is minimal and
left proper, it is primitive. In this case the sequence of sets ${\cal
P}_{n}$ defined in Section \ref{suffix} are a refining sequence of
partitions, by Corollary \ref{clopenpartition}, which generate the
topology of $(X_{\btau}, \sigma)$.  As in the proof of Proposition 16,
\cite{dhs}, we construct an isomorphism $F:X_{\btau}\rightarrow
X_{\cal B}$ by having the path $F({\by})$ passing through the vertex
in ${\cal V}_{n}$ that corresponds to the tower in ${\cal P}_{n}$
where $\by$ is located. The proof that $F$ is a conjugacy between
$(X_{\cal B}, V_{\cal B})$ and $(X_{\btau},\sigma)$ is as in the
aforementioned proof.

If $\btau$ is left proper and the branch point $\by\neq \bu$, then by Lemma \ref{inverseimagesofu}, $\btau$ (or some power of $\btau$) has to be right proper.
 Work with the substitution $\btau^{*}$ in Corollary \ref{nonfixedcorollary}. As above, $(X_{\btau^{*}},\sigma)$ is conjugate to $(X_{{\cal B}^{*}}, V_{{\cal B}^{*}})$, and since $(X_{\btau},\sigma)=(X_{\btau^{*}},\sigma)$, the result follows.\endproof

Lemma 15 in \cite{f} is used in the two sided version of this previous result, to show that (aperiodic, primitive) substitution subshifts are conjugate to a stationary adic system where there are only single edges from the vertex $v_{0}$ to any vertex in ${\cal V}_{1}$. Although the statement of this lemma  is also true for semi-proper Bratelli diagrams, it is not clear that the resulting stationary diagram corresponds to a recognizable, or injective substitution. See for example, Figure 4 in \cite{dhs}.

%with a sequence of generating partitions for
%$X_{\cal B}$. If $a\,\in {\cal V}_{n}$ and ${\bf p}$ is the minimal
%path from ${v_{0}}$ to $a$ in ${\cal B}$, then map
%$[\sigma^{k}(\btau^{n}(a))]$ to $[V_{\cal B}^{k}({\bf p})]$. If
%$[\bw]^{*}\, \in {\cal P}_{n}$ is a union of sets of the form
%$[\sigma^{k}(\btau^{n}(a))]$ let $F([\bw]^{*})$ be the corresponding
%union of paths in $X_{\cal B}$. Let ${\cal Q}_{n}:=F({\cal P}_{n})$;
%then we claim that ${\cal Q}_{n}$ is a generating sequence of
%partitions. To see this we show that if $m$ is given, then there is an $N$ such that for any $n>N$, ${\cal Q}_{n}$ is a refinement of the partition ${\cal C}_{m}$ of cylinder sets of length $m+1$.  
%Note also, that if ${\bf p} \in {\cal P}_{n}$ and ${\bf p'}\in  {\cal P}_{n'}$ with ${\bf p}\subset {\bf p}'$, then $F({\bf p})\subset F({\bf p'})$.
%Thus $F$ extends to an isomorphism $F:X_{\btau}\rightarrow X_{\cal B}$. For, if $\bx \, \in X_{\btau}$, let $F(\bx)$ be the path in $X_{\cal B}$ that passes through the vertex $V_{i}\, \in {\cal V}_{n}$ if $\bx$ lies in the 
%Check this Reem.

\section{Induced substitutions \label{inducedsubstitutions}}
We now introduce concepts needed to extend Theorem
\ref{quasiinvertibleleftproperareadic} to non-left proper
substitutions.

\subsection{Induced transformations\label{kakutani}}
Suppose that $T:X\rightarrow X$ is a minimal continuous transformation
 with $X$ a Cantor space - henceforth called a Cantor system.  If
 $U\subset X$ is a clopen set, the {\em system $(U,T_{U})$ induced by
 $(X,T)$ on $U$}, is defined by $T_{U}(\bx) = T^{n}(\bx)$, where $n$
 is the least positive natural number such that $T^{n}(x) \, \in
 U$. $T_{U} $ is well defined since $T$ is minimal, and since $X$ is
 compact, $n$ can only take a finite number of values.  The Cantor
 systems $(X_{1},T_{1})$ and $(X_{2}, T_{2})$ are (topologically) {\em
 Kakutani equivalent} (\cite{Pe} is first reference of this) if there
 exist clopen sets $U_{i}\subset X_{i}$ such that the respective
 induced systems are isomorphic. Topological Kakutani equivalence is
 much more stringent than measurable Kakutani equivalence - for
 example, all rank one transformations are measurably Kakutani
 equivalent, but it can be shown, for example, that the {\em Chacon}
 substitution $\btau(0)=0010,$ $\btau(1)=1$ is Kakutani equivalent
 only to the substitutions $\btau_{1}(0)=001^{k}0$, $\btau(1)=1$, or
 $\btau_{1}(0)=01^{k}00$, $\btau_{2}(1)$.

   Conversely, given an induced  system $(U,\sigma_{U})$, where $(X,\sigma)$ 
is minimal, for some finite $H$, we can define the {\em height function} $h:U\rightarrow \{1,\ldots H\}$ given by $h(\bx)=k$ if and only if 
$\sigma^{i}(\bx )\not\in U$ for $1\leq i\leq k-1$ but $\sigma^{k}(\bx)\in U$.
 For each $k$ in $\{1,\ldots H\}$, and $0\leq k-1$, let $U_{k}^{i}=
\sigma^{i}(h^{-1}(\{k\}))$. Since $(X,\sigma)$ is minimal,
$X=\cup_{k=1}^{H}\cup_{j=0}^{k-1}U_{k}^{j}$.

\begin{lemma}
\label{clopenlemma}
If $(X_{\btau},\sigma)$ is quasi-invertible, with branch point in the clopen set $U$, then the sets 
$\{\{U_{k}^{j}\}: 1\leq k\leq H, 0\leq j\leq k-1\}$ form a clopen partition of $X_{\btau}$.
\end{lemma}

{\bf Proof:} It is straightforward that each $U_{k}^{0}$ is clopen, by
their definition. If $({\by}_{n}) \subset U_{k}^{i}$ and ${\by}_{n}
\rightarrow {\by}$, then for each $n$ there is some ${\bf x}_{n}\,\in
U_{k}^{0}$ with $\sigma^{k}({\bf x}_{n}) = {\bf y}_{n}$. If ${\bf x}$
is a limit point of the points ${\bf x}_{n}$, then ${\bf x}\,\in
U_{k}^{0}$ and $\by=\sigma^{i}(\bx)$. Thus $U_{k}^{i}$ is closed. If the branch point in $X_{\btau}$ is in $U$, these sets are also disjoint. The result follows.
\endproof

 If we can partition $X_{\btau}$ as in Lemma \ref{clopenlemma},  we say that $(X_{\btau},\sigma)$ is a {\em primitive} of $(U,\sigma_{U})$.
If $(U,T_{U})$ has an appropriate adic representation, it is now
possible to extend this representation to $(X_{\btau},T)$, as in the proof of
Theorem 3.8, in \cite{gps}:

\begin{theorem}
\label{gpstower}
Let $(X_{\cal B}, {V}_{\cal B})$ be an adic system, where $\cal B$ is
semi proper, simple, and with M maximal elements. Suppose that the
aperiodic M-quasi-invertible Cantor system $(Y,T)$, with branch point
$\by$ has an induced system $(Y',S)$ where $\by\, \in Y'$. If $(Y',S)$
is isomorphic to $(X_{\cal B}, V_{\cal B})$, then $(Y,T)$ is
isomorphic to $(X_{\cal B'}, V_{\cal B'})$ where ${\cal B'}$ is
obtained from ${\cal B}$ by adding or removing a finite number of
edges to ${\cal E}_{1}$, and changing the ordering on the affected
vertices and edges.  \endproof\end{theorem}

\subsection{Return words}

Let $\btau$ be a (not necessarily left proper) reasonable substitution with fixed point $\bu$, and suppose that $u_{0}=a$. A {\em return word to $a$} in $\bu$ is a word
${\bw}$ such that
\begin{enumerate}
\item
$u_{0}$ is a prefix of ${\bw}$;
\item
There is no other occurrence of $u_{0}$ in ${\bw}$;
\item
${\bw}u_{0}\, \in {\cal L}_{\btau}$.
\end{enumerate}

This definition is a special case of the definition of a return word in \cite{dua}.
The set of return words ${\cal R}$ is finite, since $\bu$ is almost
periodic.  This notion is the one-sided generalization of two sided
return words in \cite{dhs}, and we use their notation here.  Ordering
${\cal R}$ according to the order of appearance of a return word in
$\bu$, we have a bijection $\psi$ from $R:=\{1,2,\ldots |{\cal R}|\}$
to ${\cal R}$. The $\btau$-fixed point $\bu$, and so every element in
$X_{\btau}$, is a concatenation of return words, with possibly a suffix of a
return word as a prefix. Extend $\psi$ by concatenation to
$\psi:R^{\mathbb N}\rightarrow [u_{0}] \subset {\cal A}^{\mathbb
N}$. Since any $\by\,\in \psi(R^{\mathbb N})$ can only be partitioned
in one way, using Property 2 of the definition of a return word,
$\psi$ is injective.

 Let ${\cal D}(\bu) \in R^{\mathbb N}$ be
the unique sequence such that $\psi({\cal D}(\bu)) = \bu$.
Let $(Y,\sigma)$ be the subshift spanned by ${\cal D}(\bu) \in R^{\mathbb N}$.
\begin{lemma}
\label{kakutanilemma}
$(Y,\sigma)$ is isomorphic to the system induced by $(X_{\btau},\sigma)$ on
$[u_{0}]$ via the map $\psi$.
\end{lemma}

{\bf Proof:} For each $\by\, \in Y$, $\psi(\sigma \by) =
\sigma^{k}(\psi (\by)$, where $k$ is the length of the word in ${\cal
R}$ corresponding to $y_{0}$. So $\psi(\sigma^{n}(\by))\, \in [u_{0}]$
for all natural $n$, thus $\psi (Y) \subset [u_{o}]$. If $\bx \, \in
[u_{0}]$, then $\sigma^{n_k}(\bu) \rightarrow \bx$ for some
$n_{k}\rightarrow \infty$. For large $k$, $\sigma^{n_k}(\bu)\in
[u_{0}]$, so $\sigma^{n_k}(\bu) = \psi (\sigma^{j_k} ({\cal D}(\bu))\,
\in \psi (Y)$ for some $j_{k}\rightarrow \infty $. Hence $\bx\, \in
\psi(Y)$. It follows that $\psi \circ \sigma = \sigma_{[u_0]}\circ
\psi$ for $y\, \in Y$.  \endproof

 The substitution $\btau$ can be used to define a substitution on $R$:
 if ${\bw}$ is a return word, then ${\bw}$ starts with $u_{0}$, and
 since $\btau(u_{0})$ starts with $u_{0}$, so does $\btau
 ({\bw})$. Also, ${\bw}u_{0}$ is a word in $\bu$, so $\btau ({\bf
 w}u_{0})= \btau({\bw})u_{0}\ldots$ is also a word in $\bu$. Hence
 $\btau({\bw})u_{0}...$ is a concatenation of unique return words, and
 a prefix of some return word, so that $\btau ({\bw})$ is a unique
 concatenation of return words. Thus if $j\, \in R$ corresponds to
 ${\bw}\,\in {\cal R}$, and $\btau ({\bw})= {\bw}_{1}{\bw}_{2}\ldots
 {\bw}_{r}$, define $\btau_{1}(j): = i_{1}i_{2}\ldots i_{r}$ where
 $i_{k} \, \in R$ corresponds to ${\bw}_{k}\, \in {\cal R}$.

\begin{lemma}
\label{towerovertau1}
Suppose $\btau$ is reasonable with generating fixed point $\bu$. Then
\begin{enumerate}
\item 
$\btau_{1}$ (or some power of $\btau_{1}$) is left proper, minimal, injective, with 
${\cal D}(\bu)$ its generating fixed point;
\item
$\btau_{1}$ is recognizable; and
\item
 If $\btau$ is $M$-quasi-invertible, then $\btau_{1}$ is
$M$-quasi-invertible. Conversely, if $\btau$ is quasi-invertible, and
$\btau_{1}$ is $M$-quasi-invertible, then $\btau$ is
$M$-quasi-invertible.
\end{enumerate}

\end{lemma}

{\bf Proof:}
\begin{enumerate}
\item
 Note that if $1\,\in R$ corresponds to ${\bw}$ in ${\cal R}$, then
${\bw}u_{0}$ is a prefix of $\bu$, and so is $\btau^{n}(u_{0})$ for
each $n$. Choose $n$ large enough so that
$|\btau^{n}(u_{0})|>|{\bw}|$, so that ${\bw}u_{0}$ is a prefix of
$\btau^{n}(u_{0})$. If $j\,\in R$, and ${\bw'}\, \in {\cal R}$
corresponds to $j$, then ${\bw'}$ begins with $u_{0}$, so that
$\btau^{n}(u_{0})$ and so ${\bw}$ is a prefix of
$\btau^{n}({\bw'})$. Hence 1 is the first letter of $\btau_{1}(j)$.

Since $\psi\circ\btau_{1}=\btau\circ\psi$, then
$\btau_{1}(D({\bu}))=D({\bu}))$. As $\btau_{1}$ is left proper, then
it is primitive. If $\btau_{1}(i)=\btau_{1}(j)$, then either $\btau
(\alpha)=\btau(\beta)$ for some $\alpha\neq \beta$, contradicting
injectivity of $\btau$, or the conditions in Theorem
\ref{mosserecognizability} are satisfied, contradicting
recognizability of $\btau$.
\item
 Suppose that $\btau_{1}$ is not recognizable. Then by Theorem
 \ref{mosserecognizability}, there exist letters $i$ and $j$ in $R$
 and some $\by \, \in Y$ with $\btau_{1}(i)$ a proper suffix of
 $\btau_{1}(j)$, and $\btau_{1}(i)\by$ and $\btau_{1}(j)\by$ appearing
 in $Y$ with the same 1-cutting of $\by$. Suppose that $\psi(i) =
 {\bw}= w_{1}\ldots w_{l}$ and $\psi(j)={\bw}'=w_{1}'\ldots
 w_{l}'$. Find the smallest $k$ in $\{0\ldots l-1\}$ such that
 $w_{l-k}=w'_{l'-k}$. Such a $k$ exists since otherwise $\bw$ is a
 suffix of $\bw'$, and since $\bw\neq \bw'$, then $w_{1}=u_{0}$ occurs
 in $w_{2}'\ldots w_{l}'$, a contradiction to Part 2 of the definition of a return word.
Now as $\btau_{1}(i)$ is a proper suffix of $\btau_{1}(j)$, there exist words $\bv$ and $\bv'$ in $R^{+}$, with $\bv$ a suffix of $\bv'$, $\bv $ a suffix of $\btau_{1}(i)$, $\bv'$ a suffix of $\btau_{1}(j)$, and words $\bp$, $\bp'$ and $\bs \, \in {\cal A}^{+}$, with 
$\psi(\bv)=\bp\btau(w_{l-k})\bs$ and $\psi(\bv')=\bp'\btau(w_{l'-k}')\bs$. If 
$|\btau(w_{l-k})|< |\btau(w_{l'-k}')|$, then $\btau(w_{l-k})$ is a proper suffix of $\btau(w_{l'-k}')$; and if $|\btau(w_{l-k})|> |\btau(w_{l'-k}')|$, then $\btau(w_{l'-k}')$ is a proper suffix of $\btau(w_{l-k})$. Finally if $|\btau(w_{l-k})|= |\btau(w_{l'-k}')|$, then $\btau(w_{l-k})=\btau(w_{l'-k}')$, which is not possible as $\btau$ is injective on letters.

Thus (without loss of generality) $\btau(w_{l-k})$ is a proper subword of $\btau(w_{l'-k}')$, and if $\by'=\btau(w_{l-k+1}\ldots w_{l})\psi(\by)$, then $\btau(w_{l-k})\by'$ and $\btau(w_{l'-k}')\by'$ appear in $X_{\btau}$ with the same 1-cutting, a contradiction to the recognizability of $\btau$.
\item
 Note that if $\bx\in X_{\btau}$ has distinct $\sigma$-preimages
$\bx_{1}, \,\bx_{2}, \ldots \bx_{M}$, then in $[u_{0}]$, there exist
$\by_{1},\, \by_{2}, \ldots \by_{M}$ with
$\by_{i}=\sigma^{-t_{i}}\bx_{i}$, where $-t_{i}$ is the first time a
preimage of $\bx_{i}$ is in $[u_{0}]$, - that this point is
well-defined follows from the quasi-invertibility and the aperiodicity
of $(X_{\btau},\sigma)$. Let $\by=\sigma_{[u_{0}]}(\bx)$ if
$\bx\not\in [u_{0}]$, and $\bx=\by$ otherwise. Then $\by$ has $M$
distinct $\sigma_{[u_{0}]}$-preimages, $\by_{1},\ldots, \by_{M}$.
Since $\bx$ is unique, so is $\by$.  Now use the isomorphism in Lemma
\ref{kakutanilemma} to transfer this information to $(Y,\sigma)$.

Conversely, if $(X_{\btau}, \sigma)$ is $N$-quasi-invertible and $\bx$
 has $N$ $\sigma$-preimages, $\{\bx_{1},\ldots \bx_{N}\}$, then if $t$
 is the first re-entry time of $\bx$ into $[u_{0}]$, $\sigma^{t}(\bx)$
 has $N$ $\sigma_{[u_0]}$-preimages $\{\sigma^{-t_1}(\bx_{1}),\ldots
 \sigma^{-t_N}\bx_{N}\}$.  Thus $\psi^{-1}(\sigma^{t}(\bx))$ has $N$
 $\sigma$-preimages, and as $(X_{\btau_1},\sigma)$ is
 $M$-quasi-invertible, $N=M$. \endproof
\end{enumerate}

\begin{theorem}
\label{quasiinvertiblearealmostadic}

Suppose that $\btau$  is  $M$-quasi-invertible and reasonable, with a generating fixed point $\bu$, and branch point $\by$.
Then 
 there
exists a semi-proper, stationary Bratteli diagram ${\cal B}$ such that $(X_{\cal B}, V_{\cal B})$ is topologically
conjugate to $(X_{\btau}, \sigma)$.

\end{theorem}

{\bf Proof:} Using the fixed point $\bu$, and working with the induced
system $([u_{0}],\sigma_{[u_{0}]})$, we use Lemma \ref{towerovertau1}
to find $\btau_{1}$ left proper, minimal, recognizable, injective and
aperiodic, such that $([u_{0}],\sigma_{[u_{0}]})$ is conjugate to
$(X_{\btau_{1}},\sigma)$.  By Part 3 of Lemma \ref{towerovertau1},
$\btau_{1}$ is $M$-quasi-invertible.  Now using Corollary
\ref{nonfixedcorollary}, there is some $\btau^{*}$ left proper,
reasonable and quasi-invertible, with its fixed point ${\bf u}^{*}$ as
branch point, such that $(X_{\btau_{1}},\sigma)= (X_{\btau^{*}},
\sigma)$ (If $\btau_{1}$'s branch point is fixed, let $\btau_{1}=
\btau^{*}$). If ${\cal B}^{*}$ is the semi-proper ordered Bratteli
diagram associated to ${\btau^{*}}$, then $(X_{{\btau}^{*}}, \sigma)$,
and so $(X_{\btau_{1}},\sigma)$, is isomorphic to $(X_{{\cal B}^{*}},
V_{{\cal B}^{*}})$. If $\by\,\in [u_{0}],$ then by Lemma
\ref{clopenlemma}, $(X_{\btau},\sigma)$ is a primitive of
$([u_{0}],\sigma_{[u_{0}]})$, and $\by$ is mapped to the minimal
element of $(X_{{\cal B}^{*}}, V_{{\cal B}^{*}})$.  Using Theorem
\ref{gpstower}, $(X_{\btau},\sigma)$ is isomorphic to $(X_{\cal B },
V_{\cal B })$ where ${\cal B}$ is obtained from ${\cal B}^{*}$ by the
addition of a finite number of edges to ${\cal E}_{1}$.

If $\by\not\in [u_{0}]$, write $(X_{\btau},\sigma)$ as a primitive over a
sufficiently small cylinder set $V$ containing $\by$, as in Lemma
\ref{clopenlemma}, and such that if a level of this partition ${\cal
P}$ intersects $[u_{0}]$, then it is contained in $[u_{0}]$: thus
${\cal P}$ has $[u_{0}]$ split up as a union of cylinder sets defined
by words of length $N$. We make $V$ small enough so that no element in
$V$ returns to $V$ before it passes through $[u_{0}]$.  Telescope the
first $M$ levels of the Bratteli diagram ${\cal B}^{*}$, so that the
elements in ${\cal P}$ which are subsets of $[u_{0}]$ are represented
by edges from $v_{0}$ to the first level. Now add levels in ${\cal
E}_{1}$ according to how the cylinder sets in $[u_{0}]$ appear in the
clopen partition of $X_{\btau}$. For example, If $C$ and $C'$ are cylinder
sets contained in $[u_{0}]$ and appearing as $C=V_{k}^{i}$ and
$C=V_{k}^{i'}$ with $i<i'$ and no set $V_{k}^{i+1},\ldots
V_{k}^{i'-1}$ intersects $[u_{0}]$, then in the telescoped Bratteli
diagram, the edge $e'$ corresponding to $C'$ is the successor of the
edge $e$ corresponding to $C$; now insert $i'-i-1$ new edges between
$e$ and $e'$.  The fact that elements in $V$ have to pass through
$[u_{0}]$ before returning to $V$ means that you can add enough levels
to the telescoped Bratteli diagram so that the new Bratteli diagram
${\cal B}$ generates an adic system $(X_{\cal B},V_{\cal B})$
isomorphic to the original $(X_{\btau},\sigma)$.

\endproof

{\bf Examples:}
\begin{enumerate}\setcounter{enumi}{8}
\item
If $\btau$ is any reasonable suffix permutative substitution, with one
fixed point $\bu$, then it is quasi-invertible, and
$(X_{\btau},\sigma)$ is a primitive over $(X_{\btau_{1}}, \sigma)$, as
$\phi(\bu) \, \in X_{\btau_{1}}$.

\item
 If $\btau(a)=aac$, $\btau(b)=bcc$, and $\btau(c)=abc$, then
$(X_{\btau},\sigma)$ is 3-quasi-invertible, with the non-invertible
element $\by$ satisfying $ c\btau(y)=y$. Since
all $\btau$-words have a proper common suffix, Corollary
\ref{nonfixedcorollary} tells us that
$(X_{\btau},\sigma)=(X_{\btau^{*}}, \sigma)$ where
$\btau^{*}(\by)=\by$. Now apply  Theorem \ref{quasiinvertibleleftproperareadic}.

\item \label{adic}

If $\btau(a)=aac$, $\btau(b)=bcc$, $\btau(c)=adbc$ and
$\btau(d)=adbd$, then since $da\not\in {\cal L}_{\btau}$ and
$cb\,\not\in {\cal L}_{\btau}$, both fixed points are
$\sigma$-invertible. Inspection of maximal proper common suffixes of
the $\btau$-substitution words of letters in any subset ${\cal A}^{*}$
of ${\cal A}$ leads to the existence of only one branch point $\by$
satisfying $c\btau(\by)=\by$. Lemma 10 cannot directly apply because
the family of $\btau$-substitution words does not have a common
suffix.

If we consider instead the induced system $([c], \sigma_{[c]})$,
$(X_{\btau}, \sigma)$ is a primitive over $([c],
\sigma_{[c]})$. Define $\btau^{*}(a)=caa$, $ \btau^{*}(b)=c$,
$\btau^{*}(c)=cadb$ and $\btau^{*}(d)=cadbdb$. One can now prove,
similarly to the proof of Lemma \ref{nonfixedlemma2a} to show that
$\btau^{*}(\by)=\by$. Thus $([c],\sigma_{[c]})$ is itself a (left
proper, quasi-invertible, reasonable) substitution subshift, and so
has a stationary adic representation $(X_{\cal B^{*}},V_{\cal
B^{*}})$. Now Theorem \ref{gpstower} applies, so that
$(X_{\btau},\sigma)$ also has a stationary adic representation. In
this example it is relatively straightforward to find the right
$\btau^{*}$; making the adic representation of $(X_{\btau},\sigma)$
more straightforward than if we had followed the proof of Theorem
\ref{quasiinvertiblearealmostadic};
in general though it is not clear how a branch point
$\by$ can be seen as the fixed point of some $\btau^{*}$ which is obtained directly from the definition of $\btau$.

\item

\label{morse}

If $\btau(0)=01$ and $\btau(1)=10$ (the {\em Morse} substitution),
then there are two $\btau$-fixed points $\bu\, \in [0]$ and
$\overline{\bu}\, \in [1]$, and these are the only branch points. The
set of return words to 0 in $\bu$ are ${\cal R}=\{ 011,01,0\}$, and
$\btau_{1}^{2}$ is defined by $ \btau_{1}^{2}(1)=123132$,
$\btau_{1}^{2}(2)=1232$ and $\btau_{1}^{2}(3)=13$, and $\btau_{1}^{2}$
has 2 branch points, the $\btau_{1}$-fixed point, and a point $\by$
satisfying $32\btau(\by)=\by$. If we take ${\cal R}$ to be the set of
return words to 1 in $\overline\bu$, then the resulting $\btau_{1}$ is
the same.  Neither $\btau$ nor $\btau_{1}$ currently has an
appropriate adic representation, so this technique fails here. It is
possible to take return words to a different set, and obtain an adic
representation, although this technique only seems to work for a
restricted family substitutions.

\item\label{chacon}
{\em Minimal rank one subshifts } are defined by substitutions on
${\cal A}=\{0,1\}$ of the form
\[\btau(0)= 0^{n_{1}}1^{m_{1}}0^{n_{2}}1^{m_{2}}\ldots
0^{n_{k-1}}1^{m_{k-1}}0^{n_{k}} \mbox{ and }\btau (1)=1,\] where
$k<\infty$, and $n_{1}$ and $n_{k}$ are positive. The latter condition ensures that the
resulting substitution subshift is minimal, with generating fixed
point $\bu:=\lim_{n\rightarrow \infty} \btau^{n}(0).$ These systems
are equivalent to rank-one systems defined by `cutting and stacking'
where there are a bounded number of cuts and spacers added, and the
same cutting-and-stacking rule is obeyed at each stage. A
comprehensive exposition of rank one systems is given in
\cite{fe}. 
Suppose $n_{i}\neq n_{j}$ for some $i,j$, so that $\btau$ is
aperiodic. Then $(X_{\btau},\sigma)$ is quasi-invertible if and only
if $m_{1}=m_{2}= \ldots = m_{k-1}=m>0$. For, the generating fixed
point $\bu$ is a branch point, and the fixed point $1^{m_i}\bu$ is another branch
point if and only 
$m_{i}<m_{j}$ for some $j$.

If $\btau(0)= 0^{n_{1}}1^{m}0^{n_{2}}1^{m}\ldots
0^{n_{k-1}}1^{m}0^{n_{k}} $, then there are two return words to 0,
$\bw_{1}=0$ and $\bw_{2}=01^{m}$. In this case $\btau_{1}$ is defined
on $\{1,2\}$ and $\btau_{i}=\bw i$ where $\bw =
1^{n_{1}-1}21^{n_{2}-1}2\ldots 1^{n_{k}-1}$. By Lemma
\ref{towerovertau1}, $\btau_{1}$ is left proper, reasonable and
2-quasi-invertible. Thus $(X_{\btau},\sigma)$ is a primitive of
$([0],\sigma_{[0]})$ and by Theorem
\ref{quasiinvertiblearealmostadic}, $\btau$ has a semi proper adic
representation by adding $m-1$ edges from $v_{0}$ to the vertex in
${\cal V}_{1}$ corresponding to the letter `1' in the Bratteli
representation of $\btau_{1}$. For $1\leq k\leq M$, the $\btau$-fixed
points $1^{k}\bu$ live on (different) levels of the tower over
$(X_{\btau},\sigma)$.

If for some $i,j$, $m_{i}\neq m_{j}$, then $\btau_{1}$ is still of the
form $\btau_{1}(i)=\bw i$, where ${\bw}= 1^{n_{1}-1}2 1^{n_{2}-1}3
\ldots 1^{n_{k-1}-1}k1^{n_{k}-1}$, so that $\btau_{1}$ is
M-quasi-invertible, with $M\geq 3$. Thus $(X_{\btau_{1}},\sigma)$ is
quasi-invertible, even though $(X_{\btau},\sigma)$ is not, and we
cannot generate an adic representation for $\btau$ from one for
$\btau_{1}$. For example, if $\btau(0)=00100110$, then the return
words are ${\cal R}=\{0,01,011\}$ and $\btau_{1}$ is 3-almost
periodic, with $\bx^{1}, \bx^{2}, $ and $\bx^{3}$ the preimages of
${\cal D} (\bu)$. In $X_{\btau}$ though, $\phi(\bx^{1})$, and $\sigma
(\phi(\bx^{2}))$ are preimages of $\bu$, and $\phi(\bx^{2})$ and
$\sigma(\phi(\bx^{3}))$ are preimages of $\sigma (\phi(\bx^{2}))$. So
an adic representation of $(X_{\btau},\sigma)$ would need 4 maximal
elements, one of which is also minimal, and another minimal element.

\item\label{orbitequivalence}
 For any primitive substitution $\btau$ whose composition matrix has a
 rational Perron-Frobenius eigenvalue, there exists a quasi-invertible
 primitive substitution $\btau'$ whose associated substitution system
 is orbit equivalent to $(X_{\btau}, \sigma)$. To see this we use the
 results in Corollary 6.7 and Theorem 6.15 in \cite{y}. First, given
 $\btau$, with Perron value $\lambda$, we find a positive integer $d$
 so that $< \mu (E): E\subset X_{\btau}\mbox{ is clopen } >= \{ \frac{n}{m}: n\,\in
 {\mathbb Z}, m \mbox{ is a factor of } d.\lambda^{n} \mbox {for some }
 n\}$, where the former is order isomorphic to the dimension group modulo the infinitesimal subgroup of $(X_{\btau}, \sigma)$. Then given $d$ and $\lambda$, both greater than 1, the author
 describes conditions on the composition matrix of ${\btau'}$ so that
 $< \mu (E): E\subset X_{\btau'}, E \mbox{ clopen }>= \{ \frac{n}{m}: n\,\in {\mathbb
 Z}, m \mbox{ is a factor of } d.\lambda^{n} \mbox { for some } n\}$.  In particular if the substitution $\btau'$ defined on $d$ letters has the substitution matrix 

\[ \left( \begin{array}{ccccc}
\lambda^{m}& \lambda^{m} & \ldots & \lambda^{m} &(\lambda^{m}-d+1)\lambda^{m} \\
\lambda^{m}&\lambda^{m}  & \ldots & (\lambda^{m}-d+1)\lambda^{m} & \lambda^{m}\\\vdots & \vdots & \vdots &\vdots &\vdots\\
(\lambda^{m}-d+1)\lambda^{m}&\lambda^{m}  &\ldots &\lambda^{m} &\lambda^{m} \end{array} \right)\] 
then $\btau'$ has the desired dimension group (modulo the infinitesimal subgroup)
 We
 claim here that $\btau'$ can be chosen to be quasi invertible. In particular if $\btau $ is defined by 
$\btau'(1)=1^{\lambda^{m}}2^{\lambda^{m}}3^{\lambda^{m}} \ldots (d-1)^{\lambda^{m}}
d^{({\lambda^{m}}-d)\lambda^{m}} d^{{\lambda^{m}}}$,
$\btau'(2)=1^{\lambda^{m}}2^{\lambda^{m}}3^{\lambda^{m}} \ldots (d-1)^{\lambda^{m}}
(d-1)^{({\lambda^{m}}-d)\lambda^{m}} d^{{\lambda^{m}}}, \ldots 
\btau'(d)=1^{\lambda^{m}}2^{\lambda^{m}}3^{\lambda^{m}} \ldots (d-1)^{\lambda^{m}}
(1)^{({\lambda^{m}}-d)\lambda^{m}} d^{{\lambda^{m}}},$
then $\btau'$ is quasi invertible. 

\item\label{notonesidedconjugate}
 Here is an example of a quasi-invertible substitution whose two-sided
subshift conjugacy class is not the same as its one-sided subshift
conjugacy class.  Define $\btau (a) = aabaa$, $\btau (b) = abcab$, and
$\btau (c) = aabac$.  The fixed point $\bu = a\ldots$ is a 3-branch
point, and $(X_{\btau}^{\mathbb N}, \sigma)$ is quasi-invertible. Let
$\bu_{L}=\ldots u_{-2}u_{-1}$ be the left infinite sequence generated
by $\{\btau^{n}(a)\}$, ${\bf v}_{L}$ the one generated by
$\{\btau^{n}(b)\}$, and ${\bf w}_{L}$ the one generated by
$\{\btau^{n}(c)\}$. Note that if $n\neq -1$, $u_{n}=w_{n}$. The words
of length 3 in ${\cal L}_{\btau}$ are ${\cal
W}_{3}=\{aab,\,\,\,aba,\,\,\,baa,\,\,\,aaa,\,\,\,abc,\,\,\, bca,\,\,\,
cab,\,\,\,bac,\,\,\,aca,\,\,\,caa\}$. Define a local rule $\phi:{\cal
W}_{3}\rightarrow \{\alpha, \beta, \gamma, d,e,f,g,h\}$ with left and
right radius one, as $\phi (aaa)=\phi(baa)=\phi (aca) = \alpha$,
$\phi(aba)=\beta$, $\phi(caa)=\gamma$, and let $\phi$ map the
remaining words in ${\cal W}_{3}$ in a one-one fashion to
$\{d,\,\,e,\,\,f,\,\,g,\,\,h\}$. Let $\Phi$ be the shift-commuting
factor mapping corresponding to $\phi$ and Let $Y=\Phi(X_{\btau}^{\mathbb Z})$.

\begin{lemma}
The map $\Phi:X_{\btau}^{\mathbb Z}\rightarrow Y$ is injective.
\end{lemma}
{\bf Proof:} Assume that for some ${\bf x}\neq {\bf y}$, $\Phi({\bf
x})= \Phi({\bf y})$. So there is some $n$ such that
$x_{n-1}x_{n}x_{n+1}\neq y_{n-1}y_{n}y_{n+1}$, and yet both words
belong to $\{baa,\,\, aaa,\,\,aca\}$. There are three cases to
consider, all similar to prove, so we look at the case when
$x_{n-1}x_{n}x_{n+1}=baa$ and $y_{n-1}y_{n}y_{n+1}=aaa$. If $baa$
appears in a sequence here, it either appears as a subword of
$\btau(a)= aabaa$, or $\btau(ba)=abcabaabaa$ or
$\btau(bc)=abcabaabac$. On the other hand $aaa$ can only appear as a
subword of $\btau(aa)=aabaaaabaa$ or $\btau(ac)=aabaaaabac$ or 
$\btau(ab)=aabaaabcab$. This
means that $x_{n-2}x_{n-1}x_{n} = aba$ and $y_{n-2}y_{n-1}y_{n}=aaa$ or $baa$. Either way 
$\Phi({\bf x})_{n}= \beta$ and $\Phi({\bf y})_{n}=\alpha$, contradicting $\Phi({\bf
x})= \Phi({\bf y})$.
\endproof

Now Theorem 4 in \cite{dhs} tells us that $(Y,\sigma)$ is either a
substitution system,or a stationary odometer. Since $(X_{\btau},\sigma)$ is not an odometer, and $\Phi$ is a conjugacy, so $(Y,\sigma)$  is generated by some substitution
$\btau'$. Hence the 2-sided subshifts generated by $\btau$ and $\btau'$ are conjugate. 

If ${\bf x}=x_{0}x_{1}\ldots$ is a one sided sequence, let $\Phi(x)= \phi(x_{0}x_{1}x_{2}) \phi(x_{1}x_{2}x_{3})\ldots$.
Note that $\btau'$ has two 2-branch points. For, \begin{eqnarray*}\Phi({\bf
u}_{L}\cdot {\bf u} ) = \Phi(\ldots u_{-2}u_{-1})\phi(u_{-2}u_{-1}u_{0})\cdot
\phi(u_{-1}u_{0}u_{1})\Phi ({\bf u}) =\\ \Phi(\ldots u_{-2}u_{-1})\phi(aaa)\cdot
\phi(aaa)\Phi ({\bf u}) = \Phi(\ldots u_{-2}u_{-1})\alpha\cdot
\alpha\Phi ({\bf u})\end{eqnarray*}
and 
 \begin{eqnarray*}\Phi({\bf
w}_{L}\cdot {\bf u} ) = \Phi(\ldots w_{-2}w_{-1})\phi(w_{-2}w_{-1}u_{0})\cdot
\phi(w_{-1}u_{0}u_{1})\Phi ({\bf u}) = \\\Phi(\ldots w_{-2}w_{-1})\phi(aca)\cdot
\phi(caa)\Phi ({\bf u}) = \Phi(\ldots w_{-2}w_{-1})\alpha\cdot
\gamma\Phi ({\bf u}),\end{eqnarray*}
so that $\Phi ({\bf u})$ is a 2-branch point; while 
\begin{eqnarray*}\Phi({\bf
v}_{L}\cdot {\bf u} ) = \Phi(\ldots v_{-2}v_{-1})\phi(v_{-2}v_{-1}u_{0})\cdot
\phi(v_{-1}u_{0}u_{1})\Phi ({\bf u}) = \\\Phi(\ldots v_{-2}v_{-1})\phi(aba)\cdot
\phi(baa)\Phi ({\bf u}) = \Phi(\ldots v_{-2}v_{-1})\beta\cdot
\alpha\Phi ({\bf u}),\end{eqnarray*}
so that $\alpha\Phi({\bf u})$ is also a 2-branch point. Thus the one sided substitution systems generated by $\btau$ and $\btau'$ cannot be topologically conjugate.

\end{enumerate}

{\footnotesize
\bibliographystyle{alpha}
\bibliography{bibliography}
}

\end{document}